\documentclass[12pt, dvips]{amsart}

\swapnumbers\newtheorem{theorem}{Theorem}[section]

\newtheorem{corollary}[theorem]{Corollary}

\newtheorem{proposition}[theorem]{Proposition}

\newtheorem{lemma}[theorem]{Lemma}

\newtheorem{question}[theorem]{Question}

\newtheorem{definition}[theorem]{Definition}

\numberwithin{equation}{section}

\usepackage{amssymb,latexsym, amsmath, amscd, bbm, graphicx, epsfig, psfrag}

\def\s{\smallskip}
\def\m{\medskip}
\def\b{\bigskip}
\def\ni{\noindent}
\def\ra{\rightarrow}
\def\ha{\hookrightarrow}

\def\peins{\operatorname{(P1)}}
\def\pzwei{\operatorname{(P2)}}
\def\pdrei{\operatorname{(P3)}}
\def\pvier{\operatorname{(P4)}}
\def\pfuenf{\operatorname{(P5)}}
\def\psechs{\operatorname{(P6)}}
\def\psieben{\operatorname{(P7)}}

\def\RR{\mathbbm{R}}

\def\ZZ{\mathbbm{Z}}

\def\aa{\alpha}
\def\bb{\beta}
\def\gg{\gamma}
\def\dd{\delta}
\def\ee{\epsilon}

\def\ff{\varphi}
\def\ll{\lambda}
\def\oo{\omega}
\def\ss{\sigma}

\def\ca{{\mathcal A}}

\def\ch{{\mathcal H}}

\def\cp{{\mathcal P}}
\def\cq{{\mathcal Q}}
\def\cR{{\mathcal R}}

\def\supp{\operatorname{supp}\;\!}
\def\Int{\operatorname{Int}}
\def\Ham{\operatorname{Ham}}

\newcommand{\proofend}{\hspace*{\fill} $\Box$\\}

\begin{document}

\title[]{On a question of Dusa Mc\,Duff}

\date{\today}

\author{Felix Schlenk}
\address{(F.\ Schlenk, School of Mathematical Sciences, Tel Aviv
  University, Ramat Aviv, Israel 69978)}
\email{schlenk@post.tau.ac.il}

\begin{abstract}  
Consider the $2n$-dimensional closed ball $B$ of radius $1$ 
in the $2n$-dimensional symplectic cylinder $Z = D \times \RR^{2n-2}$ 
over the closed disc $D$ of radius $1$.
We construct for each $\ee >0$ a Hamiltonian deformation $\ff$ of
$B$ in $Z$ of energy less than $\ee$ such that the
area of each intersection of $\ff \left( B \right)$ with the
disc $D \times \left\{ x \right\}$, $x \in \RR^{2n-2}$, is less
than $\ee$.
\end{abstract}

\maketitle

\section{Introduction}  

\ni
We endow Euclidean space $\RR^{2n}$ with the standard symplectic form
\[
\oo_0 \,=\, \sum_{i=1}^n dx_i \wedge dy_i .
\]
A $C^{\infty}$-smooth embedding $\ff$ of an open subset $U$ of
$\RR^{2n}$ into $\RR^{2n}$ is called {\it symplectic}\, if 
$\ff^* \oo_0 \,=\, \oo_0$.
An embedding of an arbitrary subset $S$ of $\RR^{2n}$ into another
subset $S'$ of $\RR^{2n}$ is called symplectic if it extends to a
symplectic embedding of a neighbourhood of $S$ into $\RR^{2n}$.
We denote by $B^{2n}(\pi r^2)$ the closed $2n$-dimensional ball of radius
$r$ and by $Z^{2n}(\pi)$ the closed $2n$-dimensional symplectic cylinder
\[
Z^{2n}(\pi) \,=\, B^2(\pi) \times \RR^{2n-2} .
\]
Gromov's celebrated Nonsqueezing Theorem \cite{G} states that there
does not exist a symplectic embedding of the ball $B^{2n}(a)$ into the
cylinder $Z^{2n}(\pi)$ if $a>\pi$.
So fix $a \in\;]0,\pi]$.
We recall that the simply connected hull $\widehat{T}$ of a subset
$T$ of $\RR^2$ is the union of its closure $\overline{T}$ and the
bounded components of $\RR^2 \setminus \overline{T}$.
We denote by $\mu$ the Lebesgue measure on $\RR^2$, and we abbreviate
$\hat{\mu}(T) = \mu \big(\widehat{T}\big)$.
It is well-known that the Nonsqueezing Theorem is equivalent to each
of the identities
\begin{eqnarray*}
\!\!\!\!\!\!\!\!\!\:\!
a &=&
\inf_\ff \, \mu \left( p \left( \ff (B^{2n}(a))  \right) \right),  \\
\!\!\!\!\!\!\!\!\!\:\!
a &=&
\inf_\ff \, \hat{\mu} \left( p \left( \ff (B^{2n}(a))  \right) \right),  
\end{eqnarray*}
where $\ff$ varies over all symplectic embeddings of $B^{2n}(a)$ into 
$Z^{2n}(\pi)$
and where $p \colon Z^{2n}(\pi) \ra B^2 (\pi)$ is the projection,
see \cite{EG} and \cite[Corollary B.10]{Diss}.
Following \cite[Section 3]{M3}
we consider sections of the image $\ff (B^{2n}(a))$ instead of its
projection, and define 
\begin{eqnarray*}    
\,\,\,\,\,\,\,\;                       
\ss (a) &=& \inf_\ff \, \sup_x \, \mu \left( p \left(
                              \ff (B^{2n}(a)) \cap D_x \right)
                          \right),\\
\,\,\,\,\,\,\,\;                     
\hat{\ss} (a) &=& \inf_\ff \, \sup_x \, \hat{\mu} \left( p \left(
                              \ff (B^{2n}(a)) \cap D_x \right) \right),
\end{eqnarray*}
where $\ff$ again varies over all symplectic embeddings of $B^{2n}(a)$ into 
$Z^{2n}(\pi)$, 
and where $D_x \subset Z^{2n} (\pi)$ denotes the disc
$D_x = B^2 (\pi) \times \{x\}$, $x \in \RR^{2n-2}$.
Clearly,
\[
\ss (a) \,\le\, \hat{\ss}(a) \,\le\, a .
\]
It is also well-known that the Nonsqueezing Theorem is equivalent
to the identity 
\begin{eqnarray}  \label{id:shutp}                     
\hat{\ss}(\pi) \,=\, \pi .
\end{eqnarray}
Indeed, the Nonsqueezing Theorem implies that for every symplectic
embedding $\ff$ of $B^{2n}(\pi)$ into $Z^{2n}(\pi)$ there exists $x
\in \RR^{2n-2}$ such that $\ff \left( B^{2n}(\pi) \right) \cap D_x$
contains the unit circle $S^1 \times \{x\}$, 
see \cite[Lemma 1.2]{LM3}.
On her search for symplectic rigidity phenomena beyond the Nonsqueezing
Theorem,  
\text{D.\ Mc\;\!Duff} therefore
asked for lower bounds of the function $\ss (a)$ and whether
$\ss (a) \ra \pi$ as $a \ra \pi$. 
It was known to \mbox{L.\ Polterovich} that $\ss (a) / a \ra 0$ as $a \ra
0$, see again \cite{M3}. 
We shall prove
\begin{theorem}\  \label{t:11}
\begin{itemize}
\item[(i)]
$\ss(a) =0$ for all $a \in \;]0,\pi]$.
\item[(ii)]
$\hat{\ss}(a) =0$ for all $a \in \;]0,\pi[$.
\end{itemize}
\end{theorem}

\ni
The symplectic embeddings in the definition of $\ss (a)$ and
$\hat{\ss} (a)$
were not further specified. 
Following a suggestion of L.\ Polterovich, 
we next ask whether the vanishing phenomenon described by
Theorem \ref{t:11} persists if we restrict ourselves to symplectic
embeddings which are close to the identity mapping in a symplectically
relevant sense.
We denote by $\ch_c (2n)$ the set of smooth functions 
$H \colon \RR^{2n} \ra \RR$ whose support is a compact subset of  
$Z^{2n}(\pi)$.
For $H \in \ch_c(2n)$ we define the Hamiltonian vector field $X_H$
through the identities
\[
\oo_0 \left( X_H(z), \cdot \right) \,=\, dH(z), \quad z \in \RR^{2n},
\]
and denote by $\phi_H$ the time-1-map of the flow generated by $X_H$. 
Moreover, we abbreviate
\begin{equation}  \label{def:Hsi}
\left\| H \right\| \,=\, \sup_{z \in \RR^{2n}} H(z) - \inf_{z \in
    \RR^{2n}} H(z) .
\end{equation}
For each $a \in \;]0,\pi]$ we define 
\begin{eqnarray*}
\ss_H (a) &=& \inf_H \bigg\{ \sup_x \, \mu \left(
p \left( \phi_H (B^{2n}(a)) \cap D_x \right) \right) + \left\| H
         \right\| \bigg\}, \\ 
\hat{\ss}_H (a) &=& \inf_H \left\{ \sup_x \, \hat{\mu} \left(
p \left( \phi_H (B^{2n}(a)) \cap D_x \right) \right) + \left\| H
         \right\| \right\}, 
\end{eqnarray*}
where $H$ varies over $\ch_c(2n)$. 
Clearly, $\ss(a) \le \ss_H(a)$ and $\hat{\ss}(a) \le \hat{\ss}_H(a)$.
In particular, $\hat{\ss}_H (\pi) = \pi$.

\begin{theorem}\  \label{t:12}
\begin{itemize}
\item[(i)]
$\ss_H(a) =0$ for all $a \in \;]0,\pi]$.
\item[(ii)]
$\hat{\ss}_H(a) =0$ for all $a \in \;]0,\pi[$.
\end{itemize}
\end{theorem}
\ni
In order to see Theorem \ref{t:12} in its right perspective we abbreviate
\[
\Ham_c \left( Z^{2n}(\pi) \right) \,=\, \left\{ \phi_H \mid H \in \ch_c
  (2n) \right\} 
\]
and define the energy $E(\phi)$ of 
$\phi \in \Ham_c \left( Z^{2n}(\pi) \right)$ by
\[
E(\phi) \,=\, \inf \left\{ \left\| H \right\| \mid \phi = \phi_H \text{ for
    some } H \in \ch_c(2n) \right\} .
\]
In the framework of Hofer geometry
the energy of a Hamiltonian diffeomorphism is its distance from the
identity mapping,
see \cite{HZ,LM3,P}.
Notice that
\begin{eqnarray*} 
\ss_H (a) &=& \inf_{\phi} \left\{ \sup_x \, \mu \left( p \left(
            \phi (B^{2n}(a)) \cap D_x \right) \right) + E(\phi) \right\}, \\ 
\hat{\ss}_H (a) &=& \inf_{\phi} \left\{ \sup_x \, \hat{\mu} \left( p \left(
            \phi (B^{2n}(a)) \cap D_x \right) \right) + E(\phi) \right\},
\end{eqnarray*}
where $\phi$ varies over $\Ham_c \left( Z^{2n}(\pi) \right)$. 
Theorem \ref{t:12} therefore says that the vanishing phenomenon
described by Theorem \ref{t:11} persists if we restrict ourselves to
Hamiltonian diffeomorphism of $Z^{2n}(\pi)$ whose Hofer distance to
the identity mapping is arbitrarily small.

\section{Results}

\ni
We start with stating a generalization of Theorem \ref{t:11}.
We denote by $\overline{\mu}$ the outer Lebesgue measure on $\RR^2$
and by $\hat{\mu}(T) = \mu \big( \widehat{T} \big)$ the Lebesgue
measure of the simply connected hull of the subset $T$ of $\RR^2$.
For each subset $S$ of the cylinder $Z^{2n}(\pi)$ we define
\begin{eqnarray*}
\ss \left( S \right) &=&  
                 \inf_\ff \, \sup_x \, \overline{\mu}  
                             \left( p \left( \ff \left( S
                 \right) \cap D_x \right) \right) , \\
\hat{\ss} \left( S \right) &=&  
                 \inf_\ff \, \sup_x \, \hat{\mu} \left( p \left(
                   \ff \left( S \right) \cap D_x \right) \right) ,
\end{eqnarray*}
where $\ff$ varies over all symplectic embeddings of $S$ into 
$Z^{2n}(\pi)$.
We abbreviate the closed cylinder $Z^{2n}(a) = B^2(a) \times \RR^{2n}$.

\begin{theorem}\  \label{t:21}
Consider a subset $S$ of $Z^{2n}(\pi)$.
\begin{itemize}
\item[(i)]
$\ss \left( S \right)=0$.

\item[(ii)]
$\hat{\ss} \left( S \right)=0$ if $S \subset Z^{2n}(a)$ for some $a<\pi$.
\end{itemize}
\end{theorem}

\ni
In view of the identity \eqref{id:shutp} we have $\hat{\ss} (S) = \pi$
whenever $S$ contains the ball $B^{2n}(\pi)$.

\begin{question}  \label{q:1}
{\rm 
Is it true that $\hat{\ss} \left( \Int B^{2n}(\pi) \right) = \pi$?
}
\end{question}

\ni
A slightly weaker version of Theorem \ref{t:21} 
has been proved in \cite{Diss} by using a symplectic folding method.
The method used here is more elementary and can also be used to prove
a generalization of Theorem \ref{t:12}.
We denote by $\ch (2n)$ the set of smooth and bounded functions 
$H \colon \RR^{2n} \ra \RR$ whose support is contained in 
$Z^{2n}(\pi)$ and whose Hamiltonian vector field $X_H$ 
generates a flow on $\RR^{2n}$. 
The time-1-map of this flow is then again denoted by $\phi_H$.
Using the notation \eqref{def:Hsi} we define
for each subset $S$ of $Z^{2n}(\pi)$, 
\begin{eqnarray*}
\ss_H (S) &=& \inf_H \left\{ \sup_x \, \overline{\mu} \left( p \left(
            \phi_H (S) \cap D_x \right) \right) + \left\| H \right\|
                    \right\}, \\
\hat{\ss}_H (S) &=& \inf_H \left\{ \sup_x \, \hat{\mu} \left( p \left(
            \phi_H (S) \cap D_x \right) \right) + \left\| H \right\| \right\},
\end{eqnarray*}
where $H$ varies over $\ch_c(2n)$ if $S$ is bounded and over $\ch
(2n)$ if $S$ is unbounded. 
In order to state the main result of this note we need yet another definition.
\begin{definition}
{\rm
A subset $S$ of $Z^{2n}(\pi)$ is {\it partially bounded}\, if
at least one of the coordinate functions $x_2, \dots, x_n, y_2, \dots,
y_n$ is bounded on $S$.
}
\end{definition}

\begin{theorem}  \label{t:22}
Consider a partially bounded subset $S$ of $Z^{2n}(\pi)$.
\begin{itemize}
\item[(i)]
$\ss_H \left( S \right)=0$.
\item[(ii)]
$\hat{\ss}_H \left( S \right)=0$ if $S \subset Z^{2n}(a)$ for some $a<\pi$.
\end{itemize}
\end{theorem}
\ni
Of course, 
\[
\ss_H \left( Z^{2n}(\pi) \right) = \hat{\ss}_H \left( Z^{2n}(\pi)
\right) =   
\ss_H \left( \Int Z^{2n}(\pi) \right) = \hat{\ss}_H \left( \Int
  Z^{2n}(\pi) \right) = \pi .
\]

\begin{question}  \label{q:2}
{\rm
Is it true that $\ss_H \left( Z^{2n}(a) \right) = \hat{\ss}_H \left(
  Z^{2n}(a) \right) =a$ for all $a \in \:]0,\pi]$? 
}
\end{question}

\ni
Theorem \ref{t:21} and Theorem \ref{t:22} are proved in the next two
sections. In Section 5 we shall reformulate these theorems in the
language of symplectic capacities.

\subsection*{Acknowledgements}
I cordially thank David Hermann and Leonid Polterovich 
for a number of helpful discussions.
I wish to thank Tel Aviv University for its hospitality
and the Swiss National Foundation for its generous support.

\section{Proof of Theorem \ref{t:21}}

\m
\ni
The main ingredient in the proof of Theorem \ref{t:21} is a special
embedding result in dimension $4$.
We shall use coordinates $z = (u,v,x,y)$ on 
$\left( \RR^4, du \wedge dv + dx \wedge dy \right)$.
We denote by $E_{(x,y)} \subset \RR^4$ the affine plane
\[
E_{(x,y)} \,=\, \RR^2 \times \{(x,y)\},
\]
and given any subset $S$ of $\RR^4$ we abbreviate
%
%
\[
\overline{\mu} \left( S \cap E_{(x,y)} \right) = 
         \overline{\mu} \left( p \left( S \cap E_{(x,y)} \right) \right) , 
\;\;
\hat{\mu} \left( S \cap E_{(x,y)} \right) = 
         \hat{\mu} \left( p \left( S \cap E_{(x,y)} \right) \right). 
\]
Fix an integer $k \ge 2$. We set
\[
\ee \,=\, \frac{\pi}{k}, \qquad 
\dd \,=\, \frac{\ee}{4k}, 
\]
and we 
define closed rectangles $P$, $P'$ and $Q$ in
$\RR^2(u,v)$ by
\begin{eqnarray*}
   \begin{array}{lcl}
\quad P  &=& [0,\pi] \times [0,1], \\ [0.3em]
\quad P' &=& [\dd, \pi-\dd] \times [\dd, 1-\dd], \\ [0.3em]
\quad Q  &=& [3\dd, \pi-3\dd] \times [3\dd, 1-3\dd].
   \end{array}
\end{eqnarray*}
We abbreviate the support of a map $\ff \colon \RR^4 \ra \RR^4$ by
\[
\supp \ff \,=\, \overline{ \left\{ z \in \RR^4 \mid \ff (z) \neq z \right\} }.
\]

\begin{proposition}  \label{p:21}
There exists a symplectomorphism $\ff$ of $\RR^4$ such that
$
\supp \ff \subset P' \times \RR^2  
$
and such that for each $(x,y) \in \RR^2$,
\begin{eqnarray}
\mu \left( \ff \left( P' \times \RR \times [0,1] \right)
    \cap E_{(x,y)} \right) &\le& 2\ee,  \label{est:sup211} \\
\hat{\mu} \left( \ff \left( Q \times \RR \times [0,1] \right)
    \cap E_{(x,y)} \right) \:&\le& 2\ee.  \label{est:sup212}
\end{eqnarray}
\end{proposition}

\proof
We define closed rectangles $R$, $R'$ and $R''$ in
$\RR^2(u,v)$ by
\begin{eqnarray*}
   \begin{array}{lcl}
\quad R &=& [0,\ee] \times [0,1], \\  [0.3em]
\quad R' &=& [\dd, \ee-\dd] \times [\dd, 1-\dd], \\  [0.3em]
\quad R'' &=& [2\dd, \ee-2\dd] \times [2\dd, 1-2\dd],
   \end{array}
\end{eqnarray*}
and we define closed rectangular annuli $A$ and $A'$ in
$\RR^2(u,v)$ by
\[
A \,=\, \overline{R \setminus R'}, \quad 
A' \,=\, \overline{R' \setminus R''} .
\]

\begin{figure}[h] 
 \begin{center}
  \psfrag{u}{$u$}
  \psfrag{v}{$v$}
  \psfrag{1}{$1$}
  \psfrag{d}{$\dd$}
  \psfrag{e}{$\ee$}
  \psfrag{A}{$A$}
  \psfrag{A'}{$A'$}
  \leavevmode\epsfbox{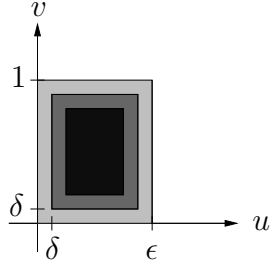}
 \end{center}
 \caption{The decomposition $R = A \cup A' \cup R''$.} 
 \label{fig1}
\end{figure}
%
%

\ni
Then $R = A \cup A' \cup R''$,
cf.\ Figure \ref{fig1}.

We choose smooth cut off functions $f_1, f_2 \colon \RR \ra [0,1]$
such that
\begin{eqnarray*}
f_1(t) &=& 
   \left\{ \begin{array}{ll}
           0,   &  t \notin [\dd, \ee - \dd], \\
           1,   &  t \in [2\dd, \ee - 2\dd],
        \end{array}
   \right. \\
f_2(t) &=& 
   \left\{ \begin{array}{ll}
           0,   &  t \notin [\dd, 1 - \dd], \\
           1,   &  t \in [2\dd, 1 - 2\dd],
        \end{array}
   \right. \\
\end{eqnarray*}
and we define the smooth function $H \colon \RR^4 \ra \RR$ by
\[
H(u,v,x,y) \,=\, - f_1(u) f_2(v) (1+\ee) x .
\]
The Hamiltonian vector field $X_H$ of $H$ is given by
\begin{equation}  \label{id:XH}
X_H(u,v,x,y) \,=\, (1+\ee) 
\left( \begin{array}{c}
- f_1(u) f_2'(v) x \\ [0.2em]
\;\;\: f_1'(u) f_2(v) x \\ [0.2em]
\;\;\: 0 \\ [0.2em]
\: f_1(u) f_2(v) 
        \end{array} 
\right) .
\end{equation}
The time-1-map $\phi_H$ has the following properties.

\m
\ni
(P1) 
$\supp \phi_H \subset R' \times \RR^2$, \\
(P2) 
$\phi_H$ fixes $A \times \RR^2$, \\
(P3) 
$\phi_H$ embeds $A' \times \RR^2$ into $A' \times \RR^2$, \\
(P4) 
$\phi_H$ translates $R'' \times \RR^2$ by $(1+\ee)1_y$.

\s
\ni
where we abbreviated $1_y = (0,0,0,1)$.

\m
For each subset $T$ of $\RR^2(u,v)$ and each $i \in \{ 1,\dots,k \}$ we
define the translate $T_i$ of $T$ by
\[
T_i \,=\, \left\{ (u+(i-1)\ee, v) \mid (u,v) \in T \right\} .
\]
With this notation we have
\[
P \,=\, \bigcup_{i=1}^k R_i \,=\, 
                           \bigcup_{i=1}^k A_i \cup A_i' \cup
                           R_i'' ,
\]

\begin{figure}[h] 
 \begin{center}
  \psfrag{u}{$u$}
  \psfrag{v}{$v$}
  \psfrag{1}{$1$}
  \psfrag{d}{$\dd$}
  \psfrag{e}{$\ee$}
  \psfrag{p}{$\pi$}
  \leavevmode\epsfbox{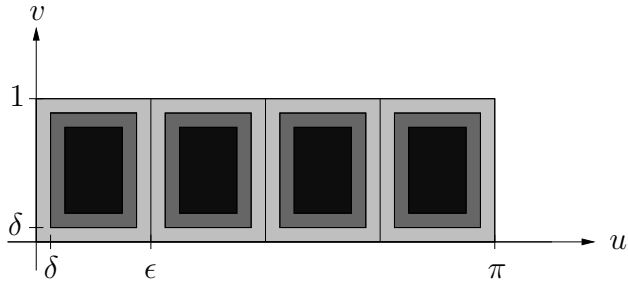}
 \end{center}
 \caption{The decomposition $P = \bigcup_{i=1}^k R_i = 
              \bigcup_{i=1}^k A_i \cup A_i' \cup R_i''$ for $k=4$.}
 \label{fig2}
\end{figure}
%
%

\ni
cf.\ Figure \ref{fig2}.
Abbreviate $H_i(u,v,x,y) = iH(u-(i-1)\ee,v,x,y)$.
We define the smooth function $\widetilde{H} \colon \RR^4 \ra \RR$ by
\[
\widetilde{H} (z) \,=\, \sum_{i=1}^k  H_i (z) 
\]
and we define the symplectomorphism $\ff$ of $\RR^4$ by 
$\ff = \phi_{\widetilde{H}}$.
In view of the identity \eqref{id:XH} we see that $\ff$ is of the form
\begin{equation}  \label{id:special}
\ff (u,v,x,y) \,=\, (u',v',x,y') ,
\end{equation}
and in view of the Properties (P1)--(P4) we find

\m
\ni
$\big( \widetilde{\text{P1}} \big)$ 
  $\supp \ff \subset P' \times \RR^2$,
  \\
$\big(\widetilde{\text{P2}}\big)$ 
  $\ff$ fixes $\bigcup_{i=1}^k A_i \times \RR^2$, \\
$\big(\widetilde{\text{P3}}\big)$ 
  $\ff$ embeds $A_i' \times \RR^2$ into $A_i' \times \RR^2$, 
  $i=1, \dots, k$, \\
$\big(\widetilde{\text{P4}}\big)$ 
  $\ff$ translates $R_i'' \times \RR^2$ by $i (1+\ee) 1_y$,
  $i=1, \dots, k$. 
\begin{figure}[h] 
 \begin{center}
  \psfrag{u}{$u$}
  \psfrag{y}{$y$}
  \psfrag{d}{$2\dd$}
  \psfrag{e}{$\ee$}
  \psfrag{p}{$\pi$}
  \psfrag{1}{$1$}
  \psfrag{1e}{$1+\ee$}
  \psfrag{2e}{$2+\ee$}
  \psfrag{3e}{$3+2\ee$}
  \psfrag{4e}{$4+3\ee$}
  \psfrag{5e}{$5+4\ee$}
  \leavevmode\epsfbox{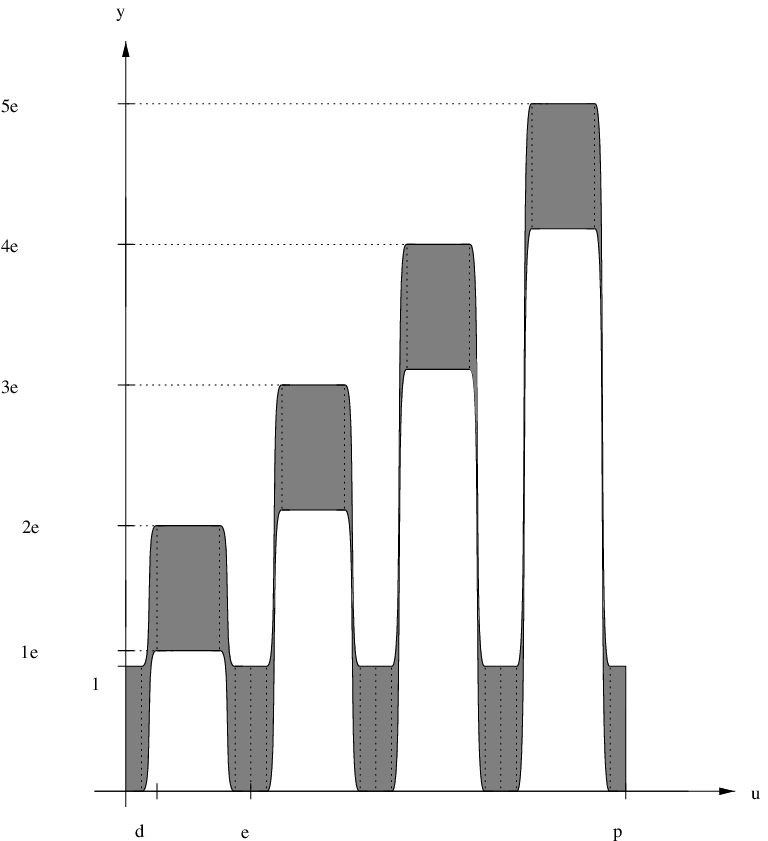}
 \end{center}
 \caption{The intersection of $\ff \left( P \times \RR 
     \times [0,1] \right)$ with a plane $\left\{
     (u,v,x,y) \mid v,x \text{ constant} \:\! \right\}$ 
                           for $v \in [2\dd, 1-2\dd]$.} 
 \label{fig3}
\end{figure}
%
%

\m
\ni
{\bf Verification of the estimates \eqref{est:sup211} and \eqref{est:sup212}} 

\s
\ni
Fix $(x, y) \in \RR^2$. 
We abbreviate 
\begin{eqnarray*}
    \begin{array}{lcl}
        \cp' &=& p \left( \ff \left( P' \times \RR \times [0,1] \right) \cap
                     E_{(x, y)} \right) ,  \\ [0.3em]
        \cq  &=& p \left( \ff \left( Q \times \RR \times [0,1] \right) \cap
                  E_{(x, y)} \right) .
    \end{array}
\end{eqnarray*}

\begin{lemma}  \label{l:disc211}
We have $\mu (\cp') \le 2 \ee$.
\end{lemma}

\proof
Using the definitions $\ee = \frac{\pi}{k}$ and $\dd = \frac{\ee}{4k}$
we estimate
\begin{equation}  \label{est:mu}
\mu \left( A_i \cup A_i' \right) \,=\, \ee - (\ee-4\dd)(1-4\dd)
\,\le\, \frac{\ee}{k}, \quad i=1, \dots, k .
\end{equation}

\m
\ni
{\bf Case A: $y \in [i^*(1+\ee), i^*(1+\ee)+1]$.}
According to Properties 
$\big( \widetilde{\text{P2}}\big)$--$\big( \widetilde{\text{P4}}\big)$  
we have $\cp' \cap R_i'' = \emptyset$ if $i \neq i^*$, and so 
\[
\cp' \,\subset\, R_{i^*} \cup \bigcup_{i=1}^k A_i \cup A_i' .
\]
Together with the estimate \eqref{est:mu} we therefore find
\begin{equation}  \label{est:p'1}
\mu \big( \cp' \big) \,\le\, \ee + k \frac{\ee}{k} \,=\, 2\ee .
\end{equation}

\s
\ni
{\bf Case B: $y \notin \bigcup_{i=1}^k [i(1+\ee), i(1+\ee)+1]$.}
According to Properties 
$\big( \widetilde{\text{P2}}\big)$--$\big( \widetilde{\text{P4}}\big)$  
we have $\cp' \cap R_i'' = \emptyset$ for all $i$, and so 
\[
\cp' \,\subset\, \bigcup_{i=1}^k A_i \cup A_i' .
\]
Therefore,
\begin{equation}  \label{est:p'2}
\mu \left( \cp' \right) \,\le\, \ee .
\end{equation}
The estimates \eqref{est:p'1} and \eqref{est:p'2} yield
that $\mu \left( \cp' \right) \le 2 \ee$.
\proofend

\begin{lemma}  \label{l:disc212}
We have $\hat{\mu} (\cq) \le 2 \ee$.
\end{lemma}

\proof
In view of the special form \eqref{id:special} of the map $\ff$ we have
\begin{equation*}  
\cq \,=\, p \big( \ff \left( Q \times \{x\} \times [0,1] \right) 
                                      \cap E_{(x, y)} \big) .
\end{equation*}
For $i = 1, \dots, k$ we abbreviate the intersections 
\begin{equation}  \label{d:AAR}
\ca_i  \,=\, Q \cap A_i, \quad 
\ca_i' \,=\, Q \cap A_i', \quad 
\cR_i'' \,=\, Q \cap R_i'' .
\end{equation}
\begin{figure}[h] 
 \begin{center}
  \psfrag{u}{$u$}
  \psfrag{v}{$v$}
  \psfrag{1}{$1$}
  \psfrag{p}{$\pi$}
  \psfrag{d}{$3\dd$}
  \leavevmode\epsfbox{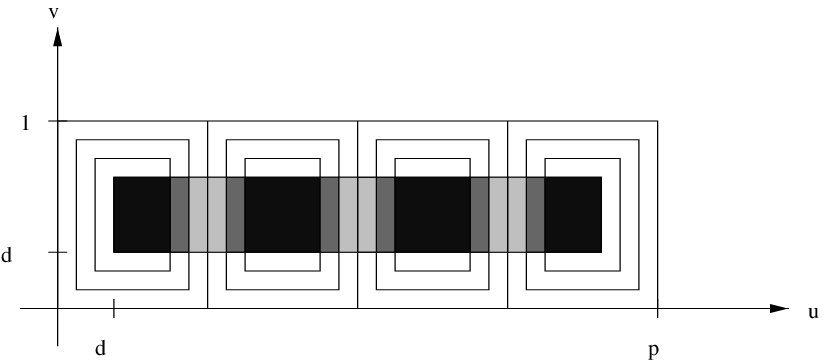}
 \end{center}
 \caption{The subsets $\ca_i$, $\ca_i'$ and $\cR_i''$ of $Q$, $i=1, \dots, 4$.} 
 \label{fig4}
\end{figure}
%
%

\ni
Each of the sets $\ca_i$ and $\ca_i'$ consists of one closed
rectangle if $i \in \{1,k\}$ and of two closed
rectangles if $i \in \{2, \dots, k-1\}$, 
cf.\ \text{Figure \ref{fig4}}.
The crucial observation in the proof is that for each $i$ the simply
connected hull of the part
\[
p \big( \ff \left( \ca_i' \times \{x\} \times [0,1] \right) 
                                      \cap E_{(x, y)} \big) 
\]
of $\cq$ is a simply connected subset of $A_i'$.
Indeed, according to property $\big( \widetilde{\text{P3}}\big)$
the closed and simply connected set
$\ff \left( \ca_i' \times \{x\} \times [0,1] \right)$ is contained
in $A_i' \times \{x\} \times \RR$, 
and so the simply connected hull of 
$\ff \left( \ca_i' \times \{ x \} \times [0,1] \right) 
\cap E_{(x,y)}$ is a simply connected subset of 
$A_i' \times \{ (x,y) \}$. 

We abbreviate by $\widehat{\cq}$ the simply connected hull of $\cq$. 

\m
\ni
{\bf Case A: $y \in [0,1]$.}
According to Properties 
$\big( \widetilde{\text{P2}}\big)$--$\big( \widetilde{\text{P4}}\big)$  
we have $\cq \cap A_i = \ca_i$ and
$\cq \cap R_i'' = \emptyset$ for all $i$.
In view of the above observation we conclude that
\[
\widehat{\cq} \,\subset\, \bigcup_{i=1}^k A_i \cup A_i' .
\]
Together with the estimate \eqref{est:mu} we therefore find
\begin{equation}  \label{est:Sa}
\mu \big( \widehat{\cq} \big) \,\le\, k \frac{\ee}{k} \,=\, \ee .
\end{equation}     

\s
\ni
{\bf Case B: $y \in [i^*(1+\ee), i^*(1+\ee)+1]$.}
According to Properties 
$\big( \widetilde{\text{P2}}\big)$--$\big( \widetilde{\text{P4}}\big)$  
we have $\cq \cap A_i = \emptyset$ for all $i$ and 
$\cq \cap R_i'' = \emptyset$ if $i \neq i^*$.
In view of the above observation we conclude that
\[
\widehat{\cq} \,\subset\, R_{i^*} \cup \bigcup_{i=1}^k A_i' .
\]
Therefore,
\begin{equation}  \label{est:Sb}
\mu \big( \widehat{\cq} \big) \,\le\, \ee + \ee \,=\, 2\ee .
\end{equation}

\s
\ni
{\bf Case C: $y \notin [0,1] \cup \bigcup_{i=1}^k [i(1+\ee), i(1+\ee)+1]$.}
According to Properties 
$\big( \widetilde{\text{P2}}\big)$--$\big( \widetilde{\text{P4}}\big)$  
we have $\cq \cap A_i = \cq \cap R_i'' = \emptyset$ for all $i$. 
In view of the above observation we conclude that
\[
\widehat{\cq} \,\subset\, \bigcup_{i=1}^k A_i' .
\]
Therefore,
\begin{equation}  \label{est:Sc}
\mu \big( \widehat{\cq} \big) \,\le\, \ee .
\end{equation}
The estimates \eqref{est:Sa}, \eqref{est:Sb} and \eqref{est:Sc} yield
that $\hat{\mu}(\cq) = \mu \big( \widehat{\cq} \big) \le 2\ee$.
This completes the proof of Lemma \ref{l:disc212}.
\proofend

In view of Lemmata \ref{l:disc211} and \ref{l:disc212} the estimates
\eqref{est:sup211} and \eqref{est:sup212} hold true.
The proof of Proposition \ref{p:21} is thus complete.
\proofend

\ni
{\bf End of the proof of Theorem \ref{t:21}\,(i)}

\s
\ni
Fix $k \ge 2$ and set $\ee = \frac{\pi}{k}$.
We choose a symplectomorphism $\aa$ of $\RR^2(u,v)$ such that $P'
\subset \aa \left( B^2(\pi) \right)$.
We refer to \cite[Lemma 2.5]{Diss} for an explicit construction.
Choose an orientation preserving diffeomorphism $f \colon \RR \ra
\;]0,1[$ and denote by $f'$ its derivative. Then the map
\[
\bb \colon \RR^2 \:\!\ra\:\! \RR \;\! \times \, ]0,1[, 
\quad
(x,y) \:\!\mapsto\:\! \left( \frac{x}{f'(y)},\, f(y)  \right)
\]
is a symplectomorphism.
We define the symplectic embedding $\Phi \colon \RR^{2n} \ha \RR^{2n}$
by
\[
\Phi \,=\, \left( (\aa^{-1} \times id) \circ \ff \circ (\aa \times
\bb) \right) \times id_{2n-4}
\]
where $\ff$ is the map guaranteed by Proposition \ref{p:21}.
Since
\begin{equation}  \label{inc:3}
\supp \ff \,\subset\, P' \times \RR^2 \,\subset\, \aa
\left(B^2(\pi)\right) \times \RR^2
\end{equation}
we have $\Phi \left( Z^{2n}(\pi) \right) \subset Z^{2n}(\pi)$.
For each subset $S$ of $Z^{2n}(\pi)$ and each point $z = (x,y,z_3,
\dots, z_n) \in \RR^{2n-2}$ we have
\begin{eqnarray*}
\Phi (S) \cap D_z 
    &\subset& \Phi \left( Z^{2n}(\pi) \right) \cap D_z \\
    &=& \left( (\aa^{-1} \times id) \circ \ff \circ (\aa
           \times \bb) \right) \left( Z^4(\pi) \right) \cap D_{(x,y)}  \\
    &\subset& \left( (\aa^{-1} \times id) \circ \ff
                  \right) \left(\aa \left( B^2(\pi) \right) \times \RR
                  \times [0,1] \right) \cap E_{(x,y)}.  
\end{eqnarray*}
Using this, the facts that $\overline{\mu}$ is monotone and $\aa^{-1}$ preserves
$\mu$, the inclusions \eqref{inc:3} and the estimates \eqref{est:sup211}
and \eqref{est:mu} we can estimate
\begin{eqnarray*}
\overline{\mu} \left( \Phi (S) \cap D_z \right) 
    &\le& \mu \left( \ff \left( \aa \left( B^2(\pi) \right) \times \RR
                  \times [0,1] \right) \cap E_{(x,y)} \right) \\
    &=& \mu \left( \ff \left( P' \times \RR \times [0,1] \right) \cap
      E_{(x,y)} \right) + \mu \left( \aa \left( B^2(\pi) \right)
      \setminus P' \right) \\
    &\le& 3\ee.
\end{eqnarray*}
Since this holds true for all $z \in \RR^{2n-2}$ 
and since $k \ge 2$ was arbitrary, we conclude that $\ss (S) =0$.

\b
\ni
{\bf End of the proof of Theorem \ref{t:21}\,(ii)}

\s
\ni
Choose $a<\pi$ so large that $S \subset Z^{2n}(a)$.
We choose $k \ge 2$ so large that $a<\mu (Q)$.
We then find a symplectomorphism $\aa$ of $\RR^2(u,v)$ such that 
\[
\aa \left( B^2(a) \right) \subset Q 
\quad \text{ and } \quad 
\aa \left( B^2(\pi) \right) \supset P' ,
\]

\begin{figure}[h] 
 \begin{center}
  \psfrag{u}{$u$}
  \psfrag{v}{$v$}
  \psfrag{1}{$1$}
  \psfrag{e}{$\ee$}
  \psfrag{p}{$\pi$}
  \psfrag{a}{$\aa$}
  \psfrag{Q}{$Q$}
  \psfrag{P}{$P'$}
  \leavevmode\epsfbox{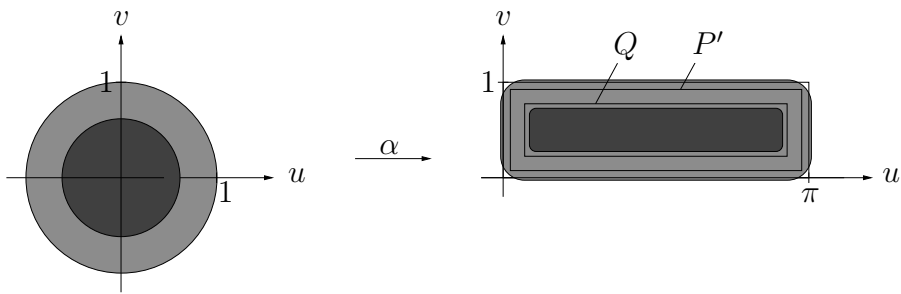}
 \end{center}
 \caption{The symplectomorphism $\aa$.} 
 \label{fig8}
\end{figure}
%
%

\ni
cf.\ Figure \ref{fig8}.
We refer again to \cite[Lemma 2.5]{Diss} for an explicit construction.
We choose a symplectomorphism $\bb \colon \RR^2 \ra \RR \:\! \times
\,]0,1[$ as above and 
define the symplectic embedding $\Phi \colon \RR^{2n} \ha \RR^{2n}$ by
\[
\Phi \,=\, \left( (\aa^{-1} \times id) \circ \ff \circ (\aa \times
\bb) \right) \times id_{2n-4} .
\]
Since
$\supp \ff \subset P' \times \RR^2 \subset \aa
\left(B^2(\pi)\right) \times \RR^2$
we have $\Phi \left( Z^{2n}(a) \right) \subset Z^{2n}(\pi)$.
For each $z = (x,y,z_3, \dots, z_n) \in \RR^{2n-2}$ we have
\begin{eqnarray*}
\Phi (S) \cap D_z 
    &\subset& \Phi \left( Z^{2n}(a) \right) \cap D_z \\
    &=& \left( (\aa^{-1} \times id) \circ \ff \circ (\aa
           \times \bb) \right) \left( Z^4(a) \right) \cap D_{(x,y)}  \\
    &\subset& \left( (\aa^{-1} \times id) \circ \ff \right) \left( Q \times \RR
                  \times [0,1] \right) \cap E_{(x,y)}.   
\end{eqnarray*}
Using this, the facts that $\hat{\mu}$ is monotone and $\aa^{-1}$ preserves
$\hat{\mu}$ and the estimate \eqref{est:sup212} we can estimate
\begin{eqnarray*}
\hat{\mu} \left( \Phi (S) \cap D_z \right) 
    &\le& \hat{\mu} \left( \ff \left( Q \times \RR \times [0,1] \right) \cap
      E_{(x,y)} \right) \\
    &\le& 2\ee.
\end{eqnarray*}
Since this holds true for all $z \in \RR^{2n-2}$ and since we can
choose $k$ as large as we like, we conclude that $\hat{\ss} (S) =0$. 
The proof of Theorem \ref{t:21} is complete.
\proofend

\section{Proof of Theorem \ref{t:22}}

\m
\ni
As in the proof of Theorem \ref{t:21} 
the main ingredient in the proof is a special
embedding result in dimension $4$.
We denote by $\ch \left( \RR^4 \right)$ the set of smooth and bounded
functions $H \colon \RR^4 \ra \RR$ whose Hamiltonian vector field $X_H$ 
generates a flow on $\RR^4$. 
The time-1-map of this flow is then again denoted by $\phi_H$, and we
abbreviate 
\begin{equation*} 
\left\| H \right\| \,=\, \sup_{z \in \RR^{2n}} H(z) - \inf_{z \in
    \RR^{2n}} H(z) .
\end{equation*}
Fix an integer $k \ge 2$ and set 
\[
\ee \,=\, \frac{\pi}{k}, \qquad 
\dd \,=\, \frac{\ee}{4k}, \qquad
\nu \,=\, \frac{\dd}{4k} .
\]
We use the notation of Section 3 and in addition 
define the closed rectangle $P^{\nu}$ in $\RR^2(u,v)$ by
\begin{equation}
P^{\nu} \,=\, [\nu, \pi-\nu] \times [\nu, 1-\nu] .
\end{equation}
We abbreviate the support of a function $H \colon \RR^4 \ra \RR$ by
\[
\supp H \,=\, \overline{ \left\{ z \in \RR^4 \mid H (z) \neq 0 \right\} }.
\]

\begin{proposition}  \label{p:22}
There exists $H \in \ch \left(\RR^4 \right)$ such that
$\supp H \subset P^{\nu} \times \RR^2$ and $\left\| H \right\| \le 2 \ee$  
and such that for each $(x,y) \in \RR^2$,
\begin{eqnarray}
\mu \left( \phi_H \left( P^{\nu} \times \RR \times [0,1] \right)
    \cap E_{(x,y)} \right) &\le& 3\ee,  \label{est:sup221} \\
\hat{\mu} \left( \phi_H \left( Q \times \RR \times [0,1] \right)
    \cap E_{(x,y)} \right) \;\,&\le& 3\ee.  \label{est:sup222}
\end{eqnarray}
\end{proposition}

\proof
As in the proof of Proposition \ref{p:21} we start with 
describing a local model of our map.
We define the closed rectangles $R$, $R'$ and $R''$ 
and the closed rectangular annuli $A$ and $A'$ as in the
proof of \text{Proposition \ref{p:21}},
and we define the closed rectangle $R^{\nu}$ in $\RR^2(u,v)$ by
\[
R^{\nu} \,=\, [\nu, \ee-\nu] \times [\nu, 1-\nu] .
\]
We also define closed intervals $I$, $I'$ and $I''$ in $\RR(x)$ by
\[
I \,=\, [0,\ee], \quad 
I' \,=\, [\dd, \ee-\dd], \quad 
I'' \,=\, [2\dd, \ee- 2\dd]
\]
and abbreviate  
\[
J  \,=\, [0,\dd] \cup [\ee-\dd, \ee], \quad
J' \,=\, [\dd, 2\dd]\: \cup \;]\ee-2\dd, \ee-\dd[.
\]
Then $I = J \cup J' \cup I''$.
We finally abbreviate 
\begin{equation}  \label{d:yy}
\check{y}_i \,=\, 1+(2i-1)\dd, \quad
\hat{y}_i \,=\, 2i - \ee + 2i \dd .
\end{equation}

\begin{lemma}  \label{l:modell}
For each $i \in \{ 1, \dots, k \}$ 
there exists a smooth  
function $H_i \colon \RR^4 \ra \RR$ with the
following properties.

\m 
\ni
$\peins_i$
$\supp H_i \subset R^{\nu} \times I \times \RR$, \\ 
$\pzwei_i$
$\phi_{H_i}$ fixes $A \times I \times [0,1]$, \\
$\pdrei_i$
$\phi_{H_i}$ embeds $A' \times I \times [0,1]$ into $A' \times I 
\times \RR$, \\
$\pvier_i$
$\phi_{H_i}$ fixes $R'' \times J \times [0,1]$, \\
$\pfuenf_i$
$\phi_{H_i}$ embeds $R'' \times J' \times [0,1]$ into
\[  
R'' \times (J \cup J') \times \RR \:\!\coprod
R'' \times I \times 
   \left( [\check{y}_i, \check{y}_i +\dd] \cup [\hat{y}_i - \ee, \hat{y}_i] \right) ,
\]
$\psechs_i$
$\phi_{H_i}$ translates $R'' \times I'' \times [0,1]$
by $2i 1_y$, \\
$\psieben_i$
$\left\| H_i \right\| \le 2 \ee$.
\end{lemma}

\begin{figure}[h] 
 \begin{center}
  \psfrag{u}{$u$}
  \psfrag{y}{$y$}
  \psfrag{1}{$1$}
  \psfrag{1e}{$1+\ee$}
  \psfrag{2}{$2$}
  \psfrag{2e}{$2+\ee$}
  \psfrag{3}{$3$}
  \psfrag{n}{$\nu$}
  \psfrag{d}{$\dd$}
  \psfrag{2d}{$2\dd$}
  \psfrag{e}{$\ee$}
  \psfrag{P}{$\phi_{H_1}$}
  \leavevmode\epsfbox{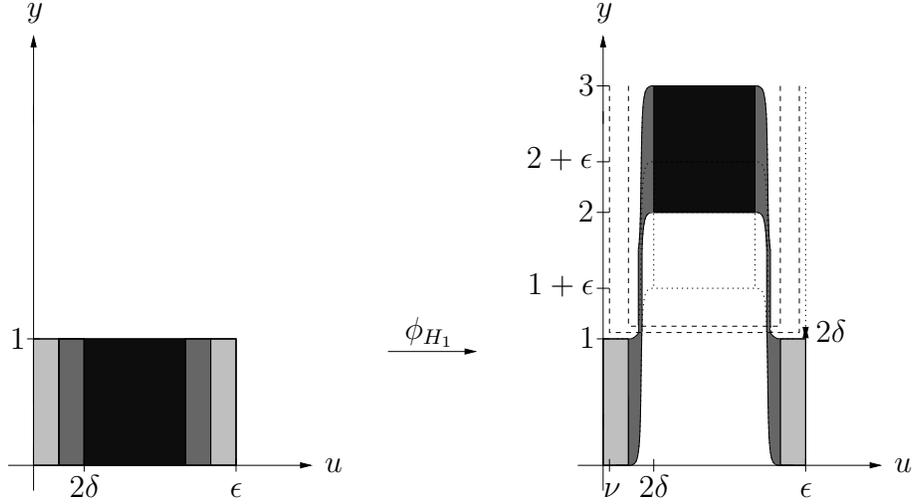}
 \end{center}
 \caption{An impression of the map $\phi_{H_1}$.} 
 \label{fig7}
\end{figure}
%
%

\proof
We shall first construct a Hamiltonian diffeomorphism $\phi_F$ of
small energy which disjoins
$R''\times I'' \times [0,1]$ from itself and shall
then construct a Hamiltonian diffeomorphism $\phi_{G_i}$ whose support
is disjoint from $R \times I \times [0,1]$ and which translates the image 
$\phi_F \left( R'' \times I'' \times [0,1] \right)$ far
along the $y$-axis.  
The composition $\phi_{G_i} \circ \phi_F \circ \phi_{G_i}^{-1}$ will be the
desired map $\phi_{H_i}$.
Both $\phi_F$ and $\phi_{G_i}$ are similar to the map $\phi_H$
constructed in the previous section, but now $F$ and $G_i$ have also
an $x$-cut off factor.
In order to make the support of $\phi_{G_i}$ disjoint from $R \times I
\times [0,1]$, the function $G_i$ must also have a $y$-cut off
factor. This will lead to technical complications.

\m
\ni
{\bf Step 1. Construction of the map $\phi_F$} 

\s
\ni
We choose smooth cut off functions $f_j \colon \RR \ra [0,1]$,
$j=1,2,3$, such that
\begin{eqnarray*}
f_1(t) &=& 
   \left\{ \begin{array}{ll}
           0,   &  t \notin [\dd,\ee-\dd], \\
           1,   &  t \in [2\dd, \ee-2\dd],
        \end{array}
   \right. \\
f_2(t) &=& 
   \left\{ \begin{array}{ll}
           0,   &  t \notin [\dd, 1 - \dd], \\
           1,   &  t \in [2\dd, 1 - 2\dd],
        \end{array}
   \right. \\
f_3(t) &=& 
   \left\{ \begin{array}{ll}
           0,   &  t \notin I', \\
           1,   &  t \in I'',
        \end{array}
   \right. 
\end{eqnarray*}
and we define the smooth function $F \colon \RR^4 \ra \RR$ by
\begin{equation*}  
F(u,v,x,y) \,=\, -f_1(u) f_2(v)f_3(x) (1+\ee) x . 
\end{equation*}
By the choice of the cut off functions $f_1$, $f_2$ and $f_3$ we have
\begin{equation}  \label{inc:XF}
\supp F\,\subset\, R' \times I' \times \RR 
\end{equation}
and since $\left| f_3(x) x \right| \le \ee-\dd$ for all $x$ we have
\begin{equation}  \label{est:2e}
\left\| F \right\| \,\le\, (1+\ee) (\ee-\dd) \,\le\, 2\ee .
\end{equation}
The Hamiltonian vector field $X_F$ of $F$ is given by
\begin{equation}  \label{id:XF}
X_F(z) \,=\, (1+\ee) 
\left( \begin{array}{c}
- f_1(u) f_2'(v) f_3(x) x \\ [0.2em]
\;\;\: f_1'(u) f_2(v) f_3(x) x \\ [0.2em]
0\\ [0.2em]
\;\;\: f_1(u) f_2(v) \left( f_3'(x) x + f_3(x) \right)
        \end{array} 
\right) .
\end{equation}
Notice that
\begin{equation}  \label{id:XFy}
X_F(z) \,=\, (1+\ee) \big( f_3'(x) x + f_3(x) \big) 1_y \quad
\text{for all } z \in R'' \times I \times \RR .
\end{equation}

\m
\ni
{\bf Step 2. Construction of the map $\phi_{G_i}$} 

\begin{figure}[h] 
 \begin{center}
  \psfrag{t}{$t$}
  \psfrag{g}{$g_4(t)$}
  \psfrag{1}{$1$}
  \psfrag{n}{$\nu$}
  \psfrag{y}{$\check{y}_i$}
  \psfrag{yd}{$\check{y}_i+\dd$}
  \leavevmode\epsfbox{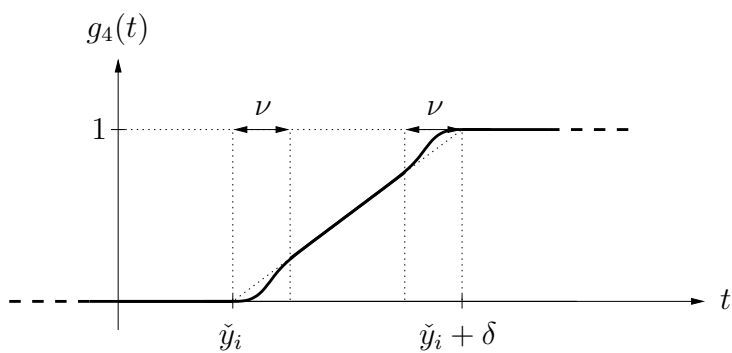}
 \end{center}
 \caption{The cut off function $g_4(t)$.} 
 \label{figg}
\end{figure}

\s
\ni
We choose smooth cut off functions $g_j \colon \RR \ra [0,1]$,
$j=1,2,3,4$, such that
\begin{eqnarray*}
g_1(t) &=& 
   \left\{ \begin{array}{ll}
           0,   &  t \notin [\nu, \ee-\nu], \\
           1,   &  t \in [\dd, \ee-\dd],
        \end{array}
   \right. \\
g_2(t) &=& 
   \left\{ \begin{array}{ll}
           0,   &  t \notin [\nu, 1-\nu], \\
           1,   &  t \in [\dd, 1 - \dd],
        \end{array}
   \right. \\
g_3(t) &=& 
   \left\{ \begin{array}{ll}
           0,   &  t \notin I, \\
           1,   &  t \in I',
        \end{array}
   \right. \\
g_4(t) &=& 
   \left\{ \begin{array}{ll}
           0,   &  t \le \check{y}_i, \\
           1,   &  t \ge \check{y}_i + \dd .
        \end{array}
   \right. 
\end{eqnarray*}
We can assume that $g_3'(t) \ge 0$ if $t \le \ee-\dd$
and that 
$g_4'(t) \ge 0$ for all $t \in \RR$ and
\begin{equation}  \label{e:g42}
g_4 (t) \,=\, \tfrac{1}{\dd} (t-\check{y}_i) \quad
  \text{if }\:\! t \in [\check{y}_i+\nu, \check{y}_i+\dd-\nu] 
\end{equation}
cf.\ Figure \ref{figg}.

%
%

\ni
We define the smooth function $G_i \colon \RR^4 \ra \RR$ by
\begin{equation*}  
G_i(u,v,x,y) \,=\, - g_1(u)g_2(v)g_3(x)g_4(y) (2i-1-\ee) x . 
\end{equation*}
The Hamiltonian vector field $X_{G_i}$ of $G_i$ is given by
\begin{equation}  \label{id:XG}
X_{G_i}(z) \,=\, (2i-1-\ee) 
\left( \begin{array}{c}
- g_1(u) g_2'(v) g_3(x) g_4(y) x \\ [0.2em]
\;\;\: g_1'(u) g_2(v) g_3(x) g_4(y) x \\ [0.2em]
- g_1(u) g_2(v) g_3(x) g_4'(y) x \\ [0.2em]
g_1(u) g_2(v) \big( g_3'(x) x + g_3(x) \big) g_4(y) 
        \end{array} 
\right) .
\end{equation}
In view of the choice of the cut off functions $g_1$ and $g_2$ we find
that for all $z \in R' \times I \times \RR$,
\begin{equation}  \label{id:XG0}
X_{G_i}(z) \,=\, (2i-1-\ee)  
\left( \begin{array}{c}
0 \\ [0.2em]
0 \\ [0.2em]
-g_3(x) g_4'(y) x \\ [0.2em]
\big( g_3'(x) x + g_3(x) \big) g_4(y)
        \end{array} 
\right) .
\end{equation}
Also notice that 
\begin{equation}  \label{inc:suppG}
\supp \phi_{G_i} \,=\, \supp \phi_{G_i}^{-1} \,\subset\, R^{\nu} \times I
\times [\check{y}_i, \infty[ . 
\end{equation}
We define the smooth function $H_i \colon \RR^4 \ra \RR$ by
\begin{equation}  \label{d:Hi}
H_i(z) \,=\, F \left( \phi_{G_i}^{-1}(z) \right) .
\end{equation}
According to the transformation law of Hamiltonian vector fields under
symplectic transformations we have
\begin{equation}  \label{e:HGFG}
\phi_{H_i} \,=\, \phi_{G_i} \circ \phi_F \circ \phi_{G_i}^{-1} .
\end{equation}

\m
\ni
{\bf Step 3. Verification of Properties $\peins_i$--$\psieben_i$} 

\s
\ni
Property $\peins_i$
follows from the inclusions \eqref{inc:XF} and \eqref{inc:suppG}.
In order to verify $\pzwei_i$--$\psieben_i$ we observe that
the inclusion \eqref{inc:suppG} and the identity \eqref{e:HGFG} imply
that
\begin{equation}  \label{id:HGF}
\phi_{H_i} (z) \,=\, \left( \phi_{G_i} \circ \phi_F \right) (z) \quad
\text{for all } z \in R \times I \times [0,1] .
\end{equation}

\s
\ni
$\pzwei_i$ and $\pvier_i$.
Assume that $z \in A \times I \times [0,1]$ or that $z \in R'' \times J
\times [0,1]$.
The inclusion \eqref{inc:XF} implies that 
$\phi_F(z) = z$.
The inclusion \eqref{inc:suppG} and the identity \eqref{id:HGF} now
imply that $\phi_{H_i}(z) =z$.

\s
\ni
$\pdrei_i$.
Assume that $z \in A' \times I \times [0,1]$.
According to the inclusion \eqref{inc:XF} and the identity
\eqref{id:XFy} we have $\phi_F(z) \in A' \times I \times \RR$. The
identities \eqref{id:XG0} and \eqref{id:HGF} now imply that
$\phi_{H_i}(z) \in A' \times I \times \RR$.

\s
\ni
$\pfuenf_i$.
Assume that $z \in R'' \times J' \times [0,1]$.
The identity \eqref{id:XFy} yields 
\[
\phi_F (z) \in R'' \times J' \times \RR .
\] 
The identity \eqref{id:XG} implies that the restriction of
$\phi_{G_i}$ to $R'' \times I \times \RR$ is of the form
\[
\phi_{G_i} (u,v,x,y) \,=\, (u,v,\ff(x,y))
\]
where $\ff$ is a symplectomorphism of $I \times \RR$.
Let $\phi_F(z) = (u_0,v_0,x_0,y_0)$. 
According to the identity \eqref{id:HGF} we need to show that
\[
\ff (x_0,y_0) \,\in\,
  (J \cup J') \times \RR \,\coprod\, 
  I \times \left( [\check{y}_i,
                  \check{y}_i+\dd] \cup [\hat{y}_i-\ee, \hat{y}_i] \right).
\]
Assume first that $y_0 \le \check{y}_i$. 
The inclusion \eqref{inc:suppG} implies that 
\[
\ff (x_0, y_0) = (x_0,y_0) \in J' \times \RR .
\]
Assume now that $y_0 \ge \check{y}_i$.
We let 
\[
\gg (t) \,=\, (x(t), y(t)),\quad t \in [0,1],
\] 
be the segment of the solution of  the system of ordinary differential
equations
\begin{eqnarray}  \label{e:sys} 
  \qquad \quad
   \left\{       
        \begin{array}{lcl}
      \dot{x}(t) &=& (2i-1-\ee) \big(\:\!\! 
                        -g_3(x(t)) g_4'( y(t)) x(t) \big) \\ [0.2em]
      \dot{y}(t) &=& (2i-1-\ee) \big(  
                        g_3'(x(t)) x(t) + g_3(x(t)) \big) g_4(y(t))     
        \end{array}
   \right. 
\end{eqnarray}
starting at $\gg(0) = (x_0, y_0)$.
Then $\gg(1) = \ff (x_0, y_0)$. 
Since $g_4'(y) \ge 0$ for all $y \in \RR$, the first equation in 
\eqref{e:sys} implies that $\dot{x}(t) \le 0$ for all $t \in [0,1]$,
and so $x(t) \le x_0 \le \ee - \dd$ for all $t \in [0,1]$.
Since $g_3'(x) \ge 0$ for all $x \le \ee - \dd$, the second equation
in \eqref{e:sys} implies that $\dot{y}(t) \ge 0$ for all $t \in
[0,1]$.

\m
\ni
{\bf Case A: $y_0 \ge \check{y}_i +\dd$.}
Since $g_4(y_0) =1$ and $\dot{y}(t) \ge 0$ we have $g_4(y(t)) =1$ for
all $t \in [0,1]$ and so $\dot{x}(t) =0$ for all $t \in [0,1]$. 
In particular, $\gg(1) \in J' \times \RR$.

\m
\ni
{\bf Case B: $y_0 \in [ \check{y}_i, \check{y}_i + \dd]$ and $x_0 \in \;]\dd,
                                                               2\dd[$.} 
Since $\dot{x}(t) \le 0$ and $\dot{y}(t) \ge 0$ for all $t \in [0,1]$,
we find that $x(1) \in [0,2\dd[$, 
and so $\gg(1) \in (J \cup J') \times \RR$.

\begin{figure}[h] 
 \begin{center}
  \psfrag{x}{$x$}
  \psfrag{y}{$y$}
  \psfrag{d}{$\dd$}
  \psfrag{e2d}{$\ee-2\dd$}
  \psfrag{ed}{$\ee-\dd$}
  \psfrag{yc}{$\check{y}_i$}
  \psfrag{ys}{$y^*$}
  \psfrag{y0}{$y_0$}
  \psfrag{yd}{$\check{y}_i+\dd$}
  \psfrag{ye}{$\hat{y}_i-\ee$}
  \psfrag{yh}{$\hat{y}_i$}
  \psfrag{C}{$C$}
  \psfrag{Ll}{$L_l$}
  \psfrag{Lr}{$L_r$}
  \psfrag{Lt}{$L_t$}
  \psfrag{Lb}{$L_b$}
  \psfrag{Lh}{$L_h$}
  \psfrag{Lv}{$L_v$}
  \psfrag{gs}{$\gg^*$}
  \psfrag{gs1}{$\gg^*(1)$}
  \psfrag{gn}{$\gg_{\nu}$}
  \psfrag{gns}{$\gg_{\nu}(t^*)$}
  \psfrag{g}{$\gg$}
  \psfrag{g1}{$\gg(1)$}
  \psfrag{gt0}{$\gg(t_0)$}
  \psfrag{n}{$\nu$}
  \leavevmode\epsfbox{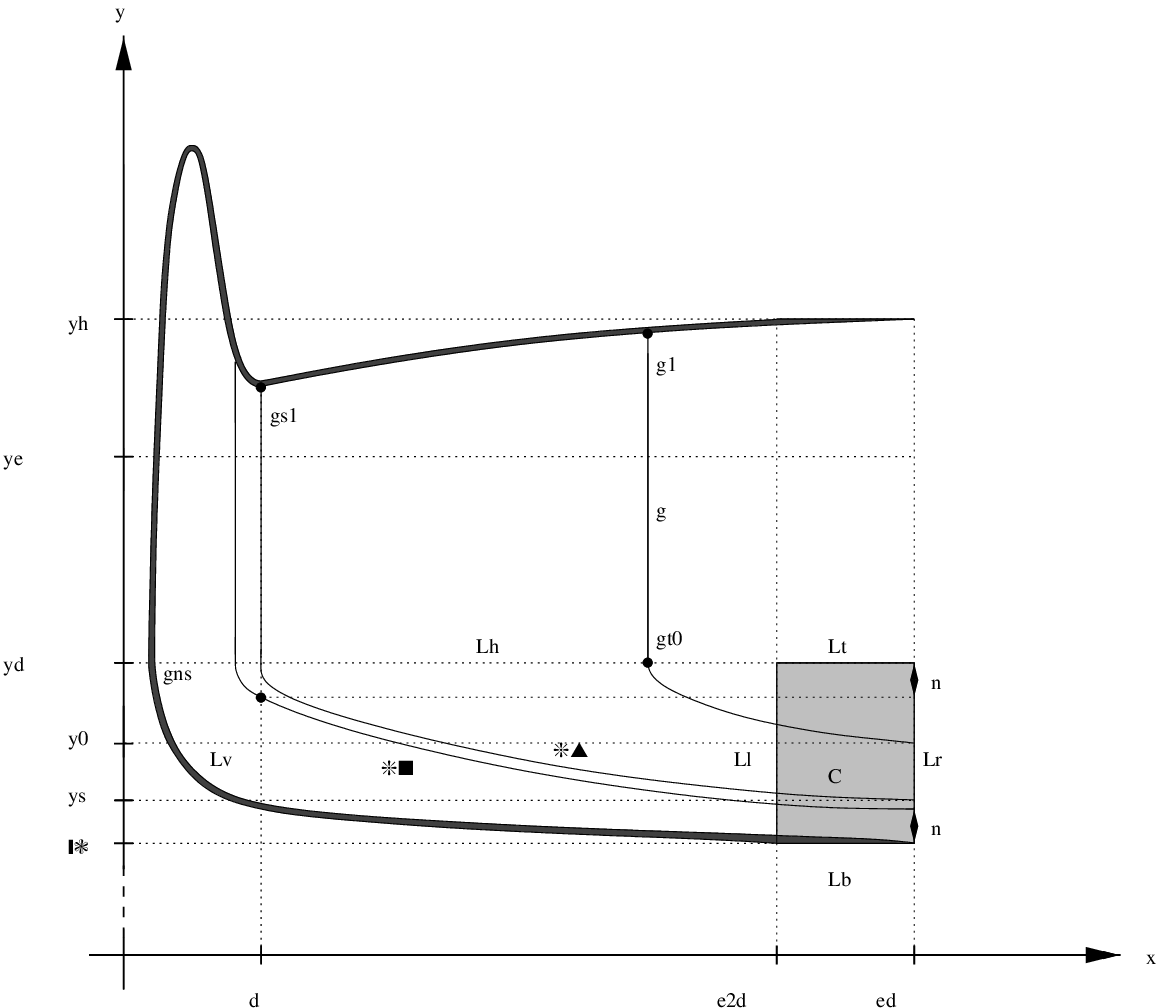}
 \end{center}
 \caption{The image $\phi(C)$, the curves $\gg_{\nu}$ and $\gg^*$
   and a curve $\gg$ starting at $(\ee-\dd, y_0)$ with $y_0>y^*$.} 
 \label{fig5}
\end{figure}

\m
\ni
{\bf Case C: $y_0 \in [ \check{y}_i, \check{y}_i+\dd]$ and $x_0 \in 
                      \;]\ee-2\dd, \ee-\dd[$.} 
We claim that
\begin{equation}  \label{inc:g1}
\gg (1) \,\in\, [0,\dd] \times \RR \,\cup\, [\dd, \ee-\dd] \times
\big( [\check{y}_i, \check{y}_i+\dd] \cup [\hat{y}_i-\ee,
\hat{y}_i] \big) .
\end{equation}
We abbreviate the closed rectangle
\[
C \,=\, [\ee - 2\dd, \ee - \dd] \times [\check{y}_i, \check{y}_i+\dd] ,
\]
and we denote the left, right, top and bottom edge of $C$ by
$L_l, \, L_r, \, L_t, \, L_b$.
It is enough to prove claim \eqref{inc:g1} for 
$(x_0,y_0) = \gg (0) \in L_l \cup L_r \cup L_t \cup L_b$, 
cf.\ Figure \ref{fig5}.
Notice that as long as 
$\gg (t) \in [\dd, \ee -\dd] \times [\check{y}_i, \check{y}_i+\dd]$, 
the system \eqref{e:sys} reads
\begin{eqnarray}  \label{deq:2} 
   \left\{ \begin{array}{lcl}
      \dot{x}(t) &=& (2i-1-\ee) \big(\:\!\! 
                        - g_4'( y(t)) x(t) \big) \\ [0.2em]
      \dot{y}(t) &=& (2i-1-\ee) g_4(y(t))     
        \end{array}
   \right. 
\end{eqnarray}
and that
\begin{eqnarray}  \label{deq:3} 
   \left\{ \begin{array}{lcl}
      \dot{x}(t) &=& 0                 \\ [0.2em]
      \dot{y}(t) &=& 2i-1-\ee     
        \end{array}
   \right. 
\end{eqnarray}
if $\gg(t) \in [\dd, \ee-\dd] \times [\check{y}_i+\dd, \infty[$.

%
%

\m
\ni
Assume $\gg(0) \in L_b$.
Then $g'_4(y_0) = g_4(y_0) =0$, and so \eqref{deq:2} implies 
\[
\gg(1) \,=\, \gg(0) \,=\, (x_0, \check{y}_i) \,\in\, [\dd, \ee-\dd] \times
[\check{y}_i, \check{y}_i+\dd] .
\]
Assume $\gg(0) \in L_t$. 
Then \eqref{deq:3} implies 
\[
\gg(1) \,=\, \gg(0) + (0, 2i-1-\ee) \,=\, (x_0, \hat{y}_i)
       \,\in\, [\dd, \ee-\dd] \times [\hat{y}_i -\ee, \hat{y}_i] .  
\] 
Assume $\gg(0) \in L_r$. 
In order to understand the locus of $\ff (L_r)$ we
abbreviate the horizontal and the vertical line 
\[
L_h \,=\, [\dd, \ee-\dd] \times \{ \check{y}_i+\dd \}
\quad \text{ and } \quad
L_v \,=\, \{ \dd \} \times [ \check{y}_i, \check{y}_i + \dd ] 
\]
and first check that the trajectory $\gg_{\nu}$ starting at $(\ee-\dd,
\check{y}_i+\nu) \in L_r$ crosses $L_v$, cf.\ Figure \ref{fig5}.
According to the choice \eqref{e:g42} the system \eqref{deq:2} reads 
\begin{eqnarray}  \label{deq:4} 
   \left\{ \begin{array}{lcl}
      \dot{x}(t) &=& (2i-1-\ee) \,\frac{1}{\dd}\, \big(\:\!\! 
                        - x(t) \big) \\ [0.2em]
      \dot{y}(t) &=& (2i-1-\ee) \,\frac{1}{\dd}\, ( y(t)
                                          -\check{y}_i )     
        \end{array}
   \right. 
\end{eqnarray}
as long as $\gg_{\nu} (t) \in [\dd, \ee-\dd] \times [ \check{y}_i+\nu,
\check{y}_i +\dd-\nu]$.
We abbreviate
\begin{equation}  \label{def:tstar}
t^* \,=\, \tfrac{\dd}{2i-1-\ee} \log \tfrac{\ee-\dd}{\dd} .
\end{equation}
Solving \eqref{deq:4} for the initial condition 
$(x_0, y_0) = (\ee-\dd, \check{y}_i + \nu)$ and using 
$\nu \frac{\ee-\dd}{\dd} = \dd -\nu$ we find that
\begin{equation}  \label{id:gn}
\gg_{\nu}(t^*) \,=\, \left( \dd, \check{y}_i + \dd -\nu \right) \,\in\, L_v .
\end{equation}
It follows that there exists a unique $y^* \in \;]\check{y}_i +\nu,
\check{y}_i+\dd[$ such that
\begin{itemize}
\item[(i)]
the curve $\gg(t)$ does not cross $L_h$ if $y_0 \in \;]\check{y}_i,
y^*[$, 
\item[(ii)]
the curve $\gg(t)$ does cross $L_h$ if $y_0 \in [y^*, \check{y}_i+\dd[$.
\end{itemize}
The trajectory $\gg^*$ starting at $(\ee-\dd, y^*)$ is shown in Figure
\ref{fig5}.

In case (i), either $\gg(t)$ does not cross the line $L_v$, in which case 
$\gg(1) \in [\dd, \ee-\dd] \times [\check{y}_i, \check{y}_i+\dd]$, or
$\gg(t)$ does cross $L_v$, in which case $\dot{x}(t) \le 0$ implies
$\gg(1) \in [0,\dd] \times \RR$. 

Assume now we are in case (ii) and that $\gg(t) = (x(t), y(t))$ is the
trajectory starting at $(y_0, \ee-\dd)$.
We define $t_0$ through the identity $\gg(t_0) \in L_h$, cf.\ Figure
\ref{fig5}. 
In order to estimate $t_0$ from above, we first notice that the
identity \eqref{id:gn} and $\check{y}_i+\nu < y_0$ imply that 
\begin{equation}  \label{est:yts}
y \left( t^* \right) \,>\, \check{y}_i + \dd - \nu .
\end{equation}
In view of the second equation in \eqref{deq:2}, the fact that $y(t)$
and $g_4(y)$ are increasing and the estimate \eqref{est:yts} we have
\[
\dot{y} (t) \,\ge\, (2i-1-\ee) \left( 1- \tfrac{1}{\dd} \nu \right)
\quad \text{ for all }\:\! t \ge t^* .
\]
Using this, the identity \eqref{def:tstar} and $\nu = \frac{\dd}{4k}$ and
$\dd = \frac{\ee}{4k}$ we can estimate
\begin{eqnarray*}
t_0 &<& t^* + \tfrac{\nu}{(2i-1-\ee)\left( 1-\tfrac{1}{\dd}\nu \right)} \\
    &<& \tfrac{1}{2i-1-\ee} \left( \dd \log (4k-1)+\dd \right) \\
    &<& \tfrac{1}{2i-1-\ee} \ee .
\end{eqnarray*}
Together with the second equation in \eqref{deq:3} and the relation
$\hat{y}_i = (\check{y}_i+\dd) + (2i-1-\ee)$ we finally find
\[
\hat{y}_i \,\ge\, y(1) \,=\, \check{y}_i + \dd + (1-t_0) (2i-1-\ee)
\,>\, \hat{y}_i - \ee 
\]
and so 
$\gg(1) \in [\dd, \ee-\dd] \times [ \hat{y}_i - \ee, \hat{y}_i]$.

\s
\ni
Assume finally $\gg(0) \in L_l$.
Since $\dot{x}(t) \le 0$ and $\dot{y}(t) \ge 0$ for all $t \in [0,1]$
we have $\gg(1) \in [0,\ee-2\dd] \times[ \check{y}_i, \infty[$.
The part of $\ff (L_l)$ contained in 
$[\dd, \ee-2\dd] \times [\check{y}_i+\dd,\infty[$ lies above the
corresponding part of $\ff (L_r)$ and below the line 
$\left\{ (x,y) \mid y = \hat{y}_i \right\}$.
Our result for $\ff (L_r)$ therefore implies that this part of $\ff
(L_l)$ is contained in $[\dd, \ee-2\dd] \times [\hat{y}_i-\ee,
\hat{y}_i]$, cf.\ \text{Figure \ref{fig5}}. We conclude that \eqref{inc:g1} 
also holds true for all points $\gg(0) \in L_l$.
This completes the verification of Property $\pfuenf_i$.

\s
\ni
$\psechs_i$.
Assume that $z \in R'' \times I'' \times [0,1]$. The identity
\eqref{id:XFy} yields 
\[
\phi_F (z) \,=\, z+(1+\ee)1_y \,\in\, R'' \times I'' \times [1+\ee, \infty[ .
\]
In view of the identities \eqref{id:HGF} and \eqref{id:XG0} and the
choice of the cut off functions $g_3$ and $g_4$ we therefore find
\[
\phi_{H_i}(z) \,=\, \phi_{G_i} \left( z+(1+\ee)1_y\right) \,=\,
z+(1+\ee)1_y + (2i-1-\ee) 1_y \,=\, z+2i 1_y .
\]

\s
\ni
$\psieben_i$.
Using the definition \eqref{d:Hi} of $H_i$ and the estimate
\eqref{est:2e} we finally estimate 
$\left\| H_i \right\| = \left\| F \right\| \le 2 \ee$. 

The proof of Lemma \ref{l:modell} is complete.
\proofend

We proceed with the proof of Proposition \ref{p:22}.
As in the previous section we define for each subset $T$ of
$\RR^2(u,v)$ and each $i \in \{1, \dots, k\}$ the translate $T_i$
of $T$ by
\[
T_i \,=\, \left\{ (u+(i-1)\ee, v) \mid (u,v) \in T \right\} ,
\]
and for each subset $X$ of $\RR(x)$ and each $i \in \{1, \dots, k\}$
and $j \in \ZZ$ we define the translate $X_{ij}$ of $X$ by
\[
X_{ij} \,=\, \left\{ x + 4(i-1)\dd + j\ee \mid x \in X \right\} .
\]
Let $H_i$ be the functions guaranteed by Lemma \ref{l:modell} and
define for each $i \in \{1,\dots,k\}$ and $j \in \ZZ$ the smooth function 
$H_{ij} \colon \RR^4 \ra \RR$ by
\[
H_{ij} (z) \,=\, H_i \left( u-(i-1)\ee, v, x-4(i-1)\dd - j\ee, y
\right) .
\]
In view of Lemma \ref{l:modell} we have

\m
\ni
$\peins_{ij}$ 
$\supp H_{ij} \subset R_i^{\nu}\times I_{ij} \times \RR$, \\
$\pzwei_{ij}$ 
$\phi_{H_{ij}}$ fixes $A_i \times I_{ij} \times [0,1]$, \\
$\pdrei_{ij}$ 
$\phi_{H_{ij}}$ embeds $A_i' \times I_{ij} \times [0,1]$ into $A_i'
                                       \times I_{ij} \times \RR$, \\
$\pvier_{ij}$ 
$\phi_{H_{ij}}$ fixes $R_i'' \times J_{ij} \times [0,1]$, \\
$\pfuenf_{ij}$ 
$\phi_{H_{ij}}$ embeds $R_i'' \times J'_{ij} \times [0,1]$ into
\[
R_i'' \times \left( J_{ij} \cup J_{ij}' \right) \times \RR \,\coprod
R_i'' \times I_{ij} \times 
   \left( [\check{y}_i, \check{y}_i+\dd] \cup [\hat{y}_i - \ee, \hat{y}_i] \right) ,
\]
$\psechs_{ij}$ 
$\phi_{H_{ij}}$ translates $R_i'' \times I''_{ij} \times [0,1]$ 
by $2i 1_y$, \\
$\psieben_{ij}$ 
$\left\| H_{ij} \right\| \le 2 \ee$.

\m
Since the sets $R^{\nu}_i \times I_{ij} \times \RR$ are
mutually disjoint, Properties $\peins_{ij}$ guarantee that the
function 
\[
H(z) \,=\, \sum_{i=1}^k \sum_{j \in \ZZ} H_{ij}(z)
\]
belongs to $\ch (4)$.
Properties $\peins_{ij}$ also imply that $\supp H \subset P^{\nu}
\times \RR^2$.
Properties $\peins_{ij}$ and $\psieben_{ij}$ imply that
\[
\left\| H \right\| \,\le\, \sup_{i,j} \left\| H_{ij}
                                     \right\| \,\le\, 2\ee .
\]

\begin{figure}[h] 
 \begin{center}
  \psfrag{u}{$u$}
  \psfrag{x}{$x$}
  \psfrag{0}{$x_0$}
  \psfrag{d}{$2\dd$}
  \psfrag{e}{$\ee$}
  \psfrag{-e}{$-\ee$}
  \psfrag{p}{$\pi$}
  \leavevmode\epsfbox{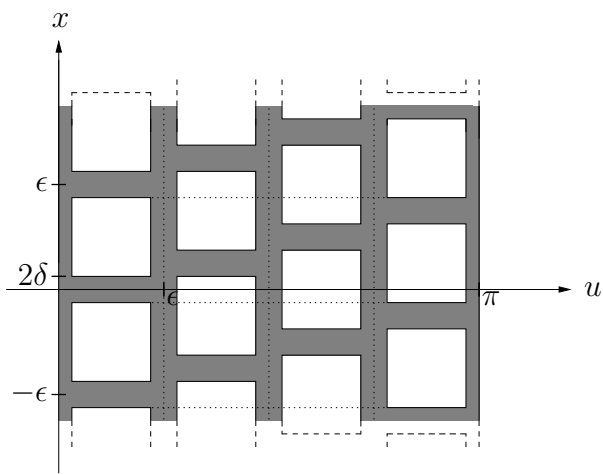}
 \end{center}
 \caption{The $(u,x)$-cut off region of the function $H$.} 
 \label{fig6}
\end{figure}
%
%


\m
\ni
{\bf Verification of the estimates \eqref{est:sup221} and \eqref{est:sup222}} 

\s
\ni
Fix $(x_0, y_0) \in \RR^2$. 
We abbreviate 
\begin{eqnarray*}
  \begin{array}{lcl}
      \cp^{\nu} &=& p \left( \phi_H \left( P^{\nu} \times \RR \times
                        [0,1] \right) \cap E_{(x_0, y_0)} \right) ,  \\ [0.3em]
      \cq       &=& p \left( \phi_H \left( Q \times \RR \times [0,1]
                        \right) \cap E_{(x_0, y_0)} \right) . 
  \end{array}
\end{eqnarray*}
Since 
\[
\RR(x) \,=\, \coprod_{i=1}^k \coprod_{j \in \ZZ} \,[-2\dd, 2\dd[_{ij}
\]
there exists a unique pair $(i_0, j_0) \in \{1, \dots, k\} \times \ZZ$
such that $x_0 \in [-2\dd, 2\dd[_{i_0j_0}$.
For $i \in \{1, \dots,k\}$ we define $j_i$ by
\begin{eqnarray*}
j_i \,=\,   
   \left\{ \begin{array}{ll}
      j_0   & \text{if }\, i \le i_0 , \\ [0.2em]
      j_0-1 & \text{if }\, i > i_0 .    
        \end{array}
   \right. 
\end{eqnarray*}
According to Properties $\peins_{ij}$ we have
\begin{eqnarray} 
\quad \cp^{\nu} \cap R_i &=& p \big( \phi_{H_{ij_i}} \left( \left(
                         P^{\nu} \cap R_i \right) \times I_{ij_i} \times
                                   [0,1] \right) \cap E_{(x_0, y_0)}
                                 \big), \label{id:Sp} \\ 
\quad \cq \cap R_i &=& p \big( \phi_{H_{ij_i}} \left( \left(
                         Q \cap R_i \right) \times I_{ij_i} \times
                                   [0,1] \right) \cap E_{(x_0, y_0)}
                                 \big),  \label{id:Sq}
\end{eqnarray}
cf.\ Figure \ref{fig6}.

\begin{lemma}  \label{l:disc221}
We have $\mu \left( \cp^{\nu} \right) \le 3 \ee$.
\end{lemma}

\proof
According to the definition \eqref{d:yy} of $\check{y}_i$ and
$\hat{y}_i$ the sets 
\[
[2i, 2i+1] \cup \left[ \check{y}_i, \check{y}_i+\dd \right] 
           \cup \left[ \hat{y}_i-\ee, \hat{y}_i \right], \;\; i=1, \dots, k ,
\]
are mutually disjoint.

\m
\ni
{\bf Case A: $y_0 \in [2i^*, 2i^*+1] \cup 
                     \left[ \check{y}_{i^*}, \check{y}_{i^*}+\dd
                     \right] \cup \left[ \hat{y}_{i^*}-\ee,
                                    \hat{y}_{i^*} \right]$.}
According to the identity \eqref{id:Sp} and 
Properties $\pzwei_{ij_i}$--$\psechs_{ij_i}$  
we have $\cp^{\nu} \cap R_i'' = \emptyset$ if $i \notin \{ i_0, i^* \}$,
and so   
\[
\cp^{\nu} \,\subset\, R_{i_0} \cup R_{i^*} \cup \bigcup_{i=1}^k A_i \cup A_i' .
\]
Together with the estimate \eqref{est:mu} we therefore find
\begin{equation}  \label{est:spn1}
\mu \left( \big( \cp^{\nu} \big) \right) \,\le\, 2\ee + k \frac{\ee}{k}
                                                           \,=\, 3\ee . 
\end{equation}     

\s
\ni
{\bf Case B: $y_0 \notin \bigcup_{i=1}^k [2i, 2i+1] \cup 
                     \left[ \check{y}_i, \check{y}_i+\dd
                     \right] \cup \left[ \hat{y}_i-\ee, \hat{y}_i \right]$.}
According to Properties $\pzwei_{ij_i}$--$\psechs_{ij_i}$  
we have 
$\cp^{\nu} \cap R_i'' = \emptyset$ if $i \neq i_0$, and so
\[
\cp^{\nu} \,\subset\, R_{i_0} \cup \bigcup_{i=1}^k A_i \cup A_i' .
\]
Therefore,
\begin{equation}  \label{est:spn2}
\mu \big( \cp^{\nu} \big) \,\le\, 2\ee .
\end{equation}
The estimates \eqref{est:spn1} and \eqref{est:spn2} 
yield that 
$\mu \left( \cp^{\nu} \right) \le 3\ee$.
\proofend

\begin{lemma}  \label{l:disc222}
We have $\hat{\mu} \left( \cq \right) \le 3 \ee$.
\end{lemma}

\proof
For $i = 1, \dots, k$ we define $\ca_i$, $\ca_i'$ and $\cR_i''$ as
in \eqref{d:AAR}.
As in the proof of Lemma \ref{l:disc212}
the crucial observation in the proof is that for each $i$ the simply
connected hull of the part
\[
p \big( \phi_{H_{ij_i}} \left( \ca_i' \times I_{ij_i} \times
                                   [0,1] \right) \cap E_{(x_0, y_0)} \big)
\]
of $\cq$ is a simply connected subset of $A_i'$.
Indeed, according to property $\pdrei_{ij_i}$
the closed and simply connected set
$\phi_{H_{ij_i}} \left( \ca_i' \times I_{ij_i} \times [0,1] \right)$ 
is contained in $A_i' \times I_{ij_i} \times \RR$, 
and so the simply connected hull of 
$\phi_{H_{ij_i}} \left( \ca_i' \times I_{ij_i} \times [0,1] \right)
\cap E_{(x_0,y_0)}$ is a simply connected subset of 
$A_i' \times \{ (x_0,y_0) \}$. 

We again abbreviate by $\widehat{\cq}$ the simply connected hull of $\cq$. 
According to the definition \eqref{d:yy} of $\check{y}_i$ and
$\hat{y}_i$ the $k+1$ sets 
\[
[0,1], \quad
[2i, 2i+1] \cup \left[ \check{y}_i, \check{y}_i+\dd \right] 
           \cup \left[ \hat{y}_i-\ee, \hat{y}_i \right], \; i=1, \dots, k ,
\]
are mutually disjoint.

\m
\ni
{\bf Case A: $y_0 \in [0,1]$.}
According to the identity \eqref{id:Sp} and 
Properties $\pzwei_{ij_i}$--$\psechs_{ij_i}$  
we have $\cq \cap A_i = \ca_i$ for all $i$ and 
$\cq \cap R_i'' = \emptyset$ if $i \neq i_0$.
In view of the above observation we conclude that
\[
\widehat{\cq} \,\subset\, R_{i_0} \cup \bigcup_{i=1}^k A_i \cup A_i' .
\]
Together with the estimate \eqref{est:mu} we therefore find
\begin{equation}  \label{est:sa}
\mu \big( \widehat{\cq} \big) \,\le\, \ee + k \frac{\ee}{k} \,=\, 2\ee .
\end{equation}     

\s
\ni
{\bf Case B: $y_0 \in [2i^*, 2i^*+1] \cup 
                     \left[ \check{y}_{i^*}, \check{y}_{i^*}+\dd
                     \right] \cup \left[ \hat{y}_{i^*}-\ee,
                                    \hat{y}_{i^*} \right]$.}
According to Properties $\pzwei_{ij_i}$--$\psechs_{ij_i}$  
we have $\cq \cap A_i = \emptyset$ for all $i$ and 
$\cq \cap R_i'' = \emptyset$ if $i \notin \{ i_0, i^* \}$.
In view of the above observation we conclude that
\[
\widehat{\cq} \,\subset\, R_{i_0} \cup R_{i^*} \cup \bigcup_{i=1}^k A_i' .
\]
Therefore,
\begin{equation}  \label{est:sb}
\mu \big( \widehat{\cq} \big) \,\le\, 2\ee + \ee \,=\, 3\ee .
\end{equation}

\s
\ni
{\bf Case C: $y_0 \notin [0,1] \cup \bigcup_{i=1}^k [2i, 2i+1] \cup 
                     \left[ \check{y}_i, \check{y}_i+\dd
                     \right] \cup \left[ \hat{y}_i-\ee, \hat{y}_i \right]$.}
According to Properties $\pzwei_{ij_i}$--$\psechs_{ij_i}$  
we have $\cq \cap A_i = \emptyset$ for all $i$ and 
$\cq \cap R_i'' = \emptyset$ if $i \neq i_0$.
In view of the above observation we conclude that
\[
\widehat{\cq} \,\subset\, R_{i_0} \cup \bigcup_{i=1}^k A_i' .
\]
Therefore,
\begin{equation}  \label{est:sc}
\mu \big( \widehat{\cq} \big) \,\le\, \ee + \ee \,=\, 2\ee .
\end{equation}
The estimates \eqref{est:sa}, \eqref{est:sb} and \eqref{est:sc} 
yield that $\hat{\mu}(\cq) = \mu \big( \widehat{\cq} \big) \le 3\ee$.
This completes the proof of Lemma \ref{l:disc222}.
\proofend

In view of Lemmata \ref{l:disc221} and \ref{l:disc222}
the estimates \eqref{est:sup221} and \eqref{est:sup222} hold true.
The proof of Proposition \ref{p:22} is thus complete.
\proofend

\ni
{\bf End of the proof of Theorem \ref{t:22}\,(i)}

\s
\ni
Consider a partially bounded subset $S$ of $Z^{2n}(\pi)$. 
There exists $i \in \{2, \dots, n\}$ and
$b>0$ such that $x_i(S) \subset [-b,b]$ or $y_i(S) \subset [-b,b]$. 
We can assume without loss of generality that $i=2$. 
If $x(S) \subset [-b,b]$, we define the symplectomorphism $\ss$ of $\RR^2(x,y)$
by $\ss (x,y) = (-y,x)$, and we let $\ss$ be the identity mapping
otherwise. Define the symplectomorphism $\tau$ of $\RR^2(x,y)$ by 
\[
\tau (x,y) \,=\, \left( 2b\,x, \frac{1}{2b}\,y + \frac{1}{2} \right) .
\]
The composition $id_2 \times (\tau \circ \ss) \times id_{2n-4}$ maps
$S$ into 
\[
B^2(\pi) \times \RR \times [0,1] \times \RR^{2n-4}.
\]
Fix $k \ge 2$.
We choose a symplectomorphism $\aa$ of $\RR^2(u,v)$ such that $P^\nu
\subset \aa \left( B^2(\pi) \right)$.
Let $H \in \ch \left( \RR^4 \right)$ be the function guaranteed by
Proposition \ref{p:22}. 
We define the smooth and bounded function $K \colon \RR^{2n} \ra \RR$
by
\begin{equation}  \label{def:K}  
K (z_1,z_2, z_3, \dots, z_n) \,=\, H \left( \aa(z_1),(\tau \circ
                              \ss)(z_2) \right) .
\end{equation}
Since 
\begin{equation}  \label{inc:HPn}  
\supp H \,\subset\, P^{\nu} \times \RR^2 \,\subset\, 
           \aa \left( B^2(\pi) \right) \times \RR^2 
\end{equation}
the support of $K$ is contained in $Z^{2n}(\pi)$, and since $\left\| H \right\|
\le 2\ee$ we have
\begin{equation}  \label{est:normK}  
\left\| K \right\| \,=\, \left\| H \right\| \,\le\, 2\ee .
\end{equation}
Moreover, the transformation law of Hamiltonian vector fields under
symplectic transformations shows that $K \in \ch (2n)$ and
\begin{equation*} 
\phi_K \,=\,  \left( (\aa \times (\tau \circ \ss))^{-1} \circ \phi_H \circ 
                  (\aa \times (\tau \circ \ss )) \right) \times id_{2n-4}.
\end{equation*}
For each subset $S$ of $Z^{2n}(\pi)$ and each point $z = (x,y,z_3,
\dots, z_n) \in \RR^{2n-2}$ we have
\begin{eqnarray*}
\phi_K (S) \cap D_z 
    &\subset& \left( (\aa \times (\tau \circ \ss ))^{-1} \circ \phi_H \right) 
           \left(\aa \left( B^2(\pi) \right) \times \RR
                  \times [0,1] \right) \cap E_{(x,y)}  \\
    &=& \left( (\aa^{-1} \times id ) \circ \phi_H \right) 
           \left(\aa \left( B^2(\pi) \right) \times \RR
                  \times [0,1] \right) \cap E_{(x',y')}
\end{eqnarray*}
where we abbreviated $(x',y') = (\tau \circ \ss) (x,y)$.
Using this, the facts that $\overline{\mu}$ is monotone and 
$\aa^{-1}$ preserves $\mu$, the inclusions \eqref{inc:HPn} 
and the estimates \eqref{est:sup221} and $\pi-\mu(P^{\nu}) \le \ee$ 
we can estimate
\begin{eqnarray*}
\overline{\mu} \left( \phi_K (S) \cap D_z \right) 
    &\le& \mu \left( \phi_H \left( \aa \left( B^2(\pi) \right) \times \RR
                  \times [0,1] \right) \cap E_{(x',y')} \right) \\
    &=& \mu \left( \phi_H \left( P^{\nu} \times \RR \times [0,1] \right) \cap
      E_{(x',y')} \right) + \mu \left( \aa \left( B^2(\pi) \right)
      \setminus P^{\nu} \right) \\
    &\le& 4\ee.
\end{eqnarray*}
Since this holds true for all $z \in \RR^{2n-2}$,
we conclude together with the estimate \eqref{est:normK} that 
\begin{equation}  \label{est:mK6}
\sup_z \overline{\mu} \left( \phi_K (S) \cap D_z \right) + \left\| K \right\|
\,\le\, 6\ee .
\end{equation}
Recall that $k \ge 2$ was arbitrary and that $\ee = \frac{\pi}{k}$. 
If $S$ is unbounded, we therefore conclude that $\ss (S) = 0$.
If $S$ is bounded, we denote by $\phi_K^t$, $t \in \RR$, the
Hamiltonian flow generated by $K$.
Since $K$ is supported in $Z^{2n}(\pi)$ and since $S$ is bounded, we find
a ball $B \subset \RR^{2n-2}$ such that
\[
\bigcup_{t \in [0,1]} \phi_K^t (S) \,\subset\, B^2(\pi) \times B .
\]
Choose a smooth compactly supported function $f \colon \RR^{2n-2} \ra
[0,1]$ such that $f |_B =1$.
The function $\widetilde{K} \colon \RR^{2n} \ra \RR$ defined by
\[
\widetilde{K} (z_1, z_2, \dots, z_n) \,=\, f(z_2, \dots, z_n) K(z_1,
\dots, z_n)
\]
belongs to $\ch_c(2n)$. Moreover, $\big\| \widetilde{K} \big\| \le
\left\| K \right\| \le 2\ee$ and $\phi_{\widetilde{K}} (S) = \phi_K (S)$. 
In view of the estimate \eqref{est:mK6} we therefore find
\[
\sup_z \overline{\mu} \left( \phi_{\widetilde{K}} (S) \cap D_z \right) + \big\|
  \widetilde{K} \big\| \,\le\, 6\ee .
\]
Since $k \ge 2$ was arbitrary, we conclude that $\ss (S) =0$. The
proof of \text{Theorem \ref{t:22}\,(i)} is complete.

\b
\ni
{\bf End of the proof of Theorem \ref{t:22}\,(ii)}

\s
\ni
Consider a partially bounded subset $S$ of $Z^{2n}(\pi)$ which is
contained in $Z^{2n}(a)$ for some $a<\pi$.
Proceeding as above we find that the composition $id_2 \times (\tau
\circ \ss) \times id_{2n-4}$ maps $S$ into
\[
B^2(a) \times \RR \times [0,1] \times \RR^{2n-4} .
\]
We choose $k \ge 2$ so large that $a<\mu (Q)$.
We then find a symplectomorphism $\aa$ of $\RR^2(u,v)$ such that
\[
\aa \left( B^2(a) \right) \subset Q 
\quad \text{ and } \quad 
\aa \left( B^2(\pi) \right) \supset P^{\nu}.
\]
We define $K \colon \RR^{2n} \ra \RR$ by formula \eqref{def:K}.
Then $K \in \ch (2n)$ and $\left\| K \right\| \le 2\ee$.
For each $z = (x,y,z_3, \dots, z_n) \in \RR^{2n-2}$ we have
\begin{eqnarray*}
\phi_K (S) \cap D_z 
    &\subset& \left( (\aa^{-1} \times id) \circ \phi_H \right) \left(
      \aa \left( B^2(a) \right) \times \RR \times [0,1] \right) 
                                         \cap E_{(x',y')}  \\ 
    &\subset& \left( (\aa^{-1} \times id) \circ \phi_H \right) 
                 \left( Q \times \RR \times [0,1] \right) \cap E_{(x',y')}.   
\end{eqnarray*}
Using this, the facts that $\hat{\mu}$ is monotone and $\aa^{-1}$ preserves
$\hat{\mu}$ and the estimate \eqref{est:sup222} we can estimate
\begin{eqnarray*}
\hat{\mu} \left( \phi_K (S) \cap D_z \right) 
    &\le& \hat{\mu} \left( \phi_H \left( Q \times \RR \times [0,1]
      \right) \cap E_{(x',y')} \right) \\
    &\le& 3\ee.
\end{eqnarray*}
Proceeding as above and recalling that we can choose $k$ as large as
we like, we conclude that $\hat{\ss} (S) =0$. 
The proof of \text{Theorem \ref{t:22}\,(ii)} is complete.
\proofend

\section{Measuring intersections by symplectic capacities}

\ni
Up to now we have measured the intersections $\ff(S) \cap D_x$ by the
outer Lebesgue measure $\overline{\mu}$ and by $\hat{\mu}$. 
There are many other ways of measuring a subset of $\RR^2$ in a
symplectic way.
We recall the

\begin{definition}
{\rm 
\cite{EH1, HZ} 
A {\it symplectic capacity on $(\RR^2, \oo_0)$} is a map
$c$ associating with each subset $T$ of $\RR^2$ 
a number $c(T) \in [0, \infty]$ in such a way that the following axioms
are satisfied.
\begin{itemize}
\item[A1.] 
{\bf Monotonicity:} $c (T) \le c(T')$ \, if there exists a
symplectomorphism $\ff$ of $\RR^2$ such that $\ff (T) \subset T'$.
\item[A2.] 
{\bf Conformality:} $c (\ll T) = \ll^2 c(T)$ \, for all  $\ll \in
                                                   \RR \setminus \{0\}$.
\item[A3.] 
{\bf Nontriviality:} $c (B^2 (\pi)) = \pi$.
\end{itemize}
A symplectic capacity $c$ on $\RR^2$ is called {\it intrinsic}\,
if it satisfies the following stronger monotonicity axiom.
\begin{itemize}
\item[A1'.] 
{\bf Monotonicity:} $c (T) \le c(T')$ \, if there exists a
symplectic embedding $\ff \colon T \ha T'$.
\end{itemize}
}
\end{definition}
\ni
Examples of intrinsic symplectic capacities on $\RR^2$ are the outer Lebesgue
measure $\overline{\mu}$, the Gromov width \cite{G} and the
Hofer--Zehnder capacity \cite{HZ}.
Examples of symplectic capacities on $\RR^2$ which are not intrinsic
are $\hat{\mu}$, the first Ekeland--Hofer capacity \cite{EH1} and the
displacement energy \cite{Ho}.
Indeed, for each of these symplectic capacities we have $c(S^1) = \pi$,
while $c(S^1) =0$ for any intrinsic symplectic capacity.
It is known that for any $a \in \;]0,\pi]$ there exists a symplectic
capacity $c$ on $\RR^2$ such that $c(S^1)=a$, see \cite[Proposition
B.11]{Diss}.
 
\s
For each subset $S$ of $Z^{2n}(\pi)$ and each symplectic capacity $c$
on $\RR^2$ we define
\[                          
\ss (S;c) \,=\, \inf_\ff \, \sup_x \, c \left( p \left(
                            \ff (S) \cap D_x \right) \right)
\]
where $\ff$ varies over all symplectic embeddings of $S$ into 
$Z^{2n}(\pi)$.
With this notation we have
$\ss (S) = \ss(S;\overline{\mu})$ and $\hat{\ss} (S) = \ss (S,\hat{\mu})$.
\begin{corollary}  \label{c:1}
Consider a subset $S$ of $Z^{2n}(\pi)$ and a symplectic capacity $c$
on $\RR^2$.
\begin{itemize}
\item[(i)]
$\ss(S;c) =0$ if $c$ is intrinsic.
\item[(ii)]
$\ss(S;c) =0$ if $S \subset Z^{2n}(a)$ for some $a <\pi$.
\end{itemize}
\end{corollary}

\proof
Consider a bounded subset $T$ of $\RR^2$.
According to \cite[Theorem B.7]{Diss} we have $c(T) \le
\overline{\mu}(T)$ for every intrinsic symplectic capacity $c$ on
$\RR^2$ and $c(T) \le \hat{\mu} (T)$ for every symplectic capacity $c$
on $\RR^2$.
Corollary \ref{c:1} thus follows from Theorem \ref{t:21}.
\proofend

For each subset $S$ of $Z^{2n}(\pi)$ and each symplectic capacity $c$
on $\RR^2$ we define
\[                          
\ss_H (S;c) \,=\, \inf_H \left\{ \sup_x \, c \left(
             p \left( \phi_H (S) \cap D_x \right) \right) + \left\| H
             \right\| \right\}
\]
where $H$ varies over $\ch_c(2n)$ if $S$ is bounded and over $\ch
(2n)$ if $S$ is unbounded. 
With this notation we have
$\ss_H (S) = \ss_H(S;\overline{\mu})$ and $\hat{\ss}_H (S) = \ss_H (S,\hat{\mu})$.
\begin{corollary}  \label{c:2}
Consider a partially bounded subset $S$ of $Z^{2n}(\pi)$ and a
symplectic capacity $c$ on $\RR^2$.
\begin{itemize}
\item[(i)]
$\ss_H(S;c) =0$ if $c$ is intrinsic.
\item[(ii)]
$\ss_H(S;c) =0$ if $S \subset Z^{2n}(a)$ for some $a <\pi$.
\end{itemize}
\end{corollary}

\proof
Corollary \ref{c:2} follows from Theorem \ref{t:22} in the same way as
Corollary \ref{c:1} followed from Theorem \ref{t:21}.
\proofend

\enddocument